\documentclass[twoside,11pt]{article}

\usepackage[abbrvbib, preprint]{jmlr2e}
\usepackage{amsmath,mathtools}
\usepackage[title]{appendix}
\usepackage[ruled]{algorithm2e}
\usepackage{matlab-prettifier}
\usepackage{multirow}
\usepackage{booktabs}
\usepackage{longtable}
\usepackage{bbm}
\usepackage{bm}
\usepackage{cite}
\usepackage[usestackEOL]{stackengine}
\usepackage{blindtext}

\newcommand{\abs}[1]{\left\lvert #1 \right\rvert}
\newcommand{\argmin}[1]{\underset{#1}{\arg\min}\,}
\newcommand{\inner}[1]{\left\langle #1 \right\rangle }
\newcommand{\norm}[1]{ \left\| #1 \right\|}

\usepackage{lastpage}

\firstpageno{1}

\begin{document}

\title{On Regularized Square-root Regression Problems:  Distributionally Robust Interpretation and Fast Computations}

\author{\name Hong T.M. Chu \email hongtmchu@u.nus.edu \\
       \addr Department of Mathematics\\
       National University of Singapore\\
       Singapore 119076
       \AND
       \name Kim-Chuan Toh \email mattohkc@nus.edu.sg \\
       \addr Department of Mathematics, and Institute of Operations Research and Analytics\\
       National University of Singapore\\
       Singapore 119076
       \AND
       \name Yangjing Zhang \email yangjing.zhang@amss.ac.cn \\
       \addr Institute of Applied Mathematics, Academy of Mathematics and Systems Science\\
       Chinese Academy of Sciences\\
       People's Republic of China 100190
       }


\maketitle

\begin{abstract}
Square-root (loss) regularized models have recently become popular in linear regression due to their nice statistical properties. Moreover, some of these models can be interpreted as the distributionally robust optimization  counterparts of the traditional least-squares regularized models. In this paper, we give a unified proof to show that any square-root regularized model whose penalty function being the sum of a simple norm and a seminorm can be interpreted as the distributionally robust optimization (DRO) formulation of the corresponding least-squares problem. In particular, the optimal transport cost in the DRO formulation is given by a certain dual form of the penalty. To solve the resulting square-root regularized model whose loss function and penalty function are both nonsmooth, we design a proximal point dual semismooth Newton algorithm and demonstrate its efficiency when the penalty is the  sparse group Lasso penalty or the  fused Lasso penalty. Extensive experiments demonstrate that our algorithm is highly efficient for solving the square-root sparse group Lasso problems and the square-root fused Lasso problems.
\end{abstract}%


\section{Introduction}

Variable selection in high dimensional feature space has played a pivotal role in contemporary statistical and machine learning studies. Let $(X_1,Y_1),\dots,(X_N,Y_N)$ be data generated from a linear regression model
\begin{equation}\label{linear-reg}
Y_i = X_i^T \beta_0 + \sigma \epsilon_i,
\end{equation}
where $X_i = (X_{i1},\dots,X_{in})^T \in \mathbb{R}^n$ is the $i$th predictor vector, $Y_i\in\mathbb{R}$ is the corresponding  response variable, and $\epsilon_i$'s are independent and identically distributed (i.i.d.) noises such that ${\rm E}(\epsilon_i) = 0$ and $ {\rm E}(\epsilon_i^2) = 1$. The vector  $\beta_0 \in \mathbb{R}^n$ is the unknown true regression vector, and $\sigma > 0$ is the unknown noise level. We let $X:= (X_1,\dots,X_N)^T \in \mathbb{R}^{N \times n}$ be the matrix of features and $Y := (Y_1,\dots,Y_N)^T \in \mathbb{R}^N$ be the response vector. In high dimensional regression problems, the
 dimension of predictors $n$ is usually much larger than the sample size $N$. With such a large number of predictors, one often prefers to select a smaller subset that is significant and highly correlated to the response for the ease of interpreting the model, as well as to avoid the issue of overfitting.
One of the most popular methods for variable selection is the classic Lasso model \citep{tibshirani1996regression}, which minimizes the residual sum of squared errors plus an $\ell_1$ norm penalty term. Further imposing problem-specific assumptions on the structures of the variables,
various variants of the classic Lasso model have been proposed. These variants are referred to as the Lasso-type models, and they include for example, the  (sparse) group Lasso model \citep{yuan2006model,friedman2010note}, the fused Lasso model \citep{tibshirani2005sparsity}, and the clustered Lasso model \citep{petry2011pairwise,she2010sparse}, to mention only a few.
To date, large amount of research have been conducted on Lasso-type models, and they have been widely applied in real applications; see for example \citep{xu2010simple,muthukrishnan2016lasso,yang2010online,jacob2009group,angelosante2009rls,rao2015classification,bazerque2011group}, just to name a few.
Most Lasso-type models incorporate different regularizers for achieving different underlying sparsity structures in the regression vector
and they solve the following generic squared-loss convex optimization problem:
\begin{equation}\label{eq:lassotype}
  \min_{\beta\in\mathbb{R}^n} \quad \left\{ \|Y - X\beta\|^2 + \lambda p(\beta) \right\},
\end{equation}
where $\|\cdot\|$ is the Euclidean norm, $p$ is the regularizer, and $\lambda > 0$ is the tuning parameter.

The theoretically optimal values of the tuning parameter $\lambda$ in the Lasso and group Lasso models  have been well estimated in \citep{bickel2009simultaneous,lounici2011oracle,meinshausen2009lasso,zhang2008sparsity}. However, an outstanding problem is that  the choice of $\lambda$ will depend on the unknown noise level $\sigma$, which is typically nontrivial to estimate. An alternative for choosing a suitable $\lambda$ for good generalization performance is by cross validation, but this procedure can be time-consuming in practice. To eliminate the need to know or estimate $\sigma$, \citet{belloni2011square}  introduced the square-root Lasso model, and subsequently \citet{bunea2013group}  extended it to the square-root group Lasso model. The  square-root (loss) regularized model takes the square-root of the residual sum of  squared errors as the loss function and any sparsity inducing norm as the regularizer. From now on we refer to the following optimization problem as the  square-root  regularized model and its solution as the  square-root  regularized estimator
\begin{equation}\label{eq:square-root}
\min_{\beta\in\mathbb{R}^n} \quad \left\{ \|Y - X\beta\| + \lambda p(\beta) \right\}.
\end{equation}
An important general step forward has been made by \citet{stucky2017sharp}, who proved the nice statistical property that the theoretically optimal $\lambda$ for an estimator given by \eqref{eq:square-root} will be independent of the unknown noise level $\sigma$ for any penalty $p$
that is a
weakly decomposable norm \citep[Definition~4.1]{van2014weakly}. Such an estimator is quite general in that $p$ can be an arbitrary weakly decomposable norm. In fact, it was stated in \citep{stucky2017sharp} that the nice statistical property is applicable for
the square-root Lasso, the square-root (sparse) group Lasso, and the square-root SLOPE \citep{bogdan2015slope} estimators. In addition, it has been also shown by \citet{MR4159660} that the choice of $\lambda$ in the square-root fused Lasso model is independent of $\sigma$. Thus far, we can see that the  square-root regularized model \eqref{eq:square-root} is
more preferable over the Lasso-type model \eqref{eq:lassotype} in terms of their statistical properties.

Recently, there has been renewed interest in the  square-root regularized  model \eqref{eq:square-root} and its interpretation from the perspective of
distributionally robust optimization (DRO). This connection gives a new probabilistic explanation of the penalty level $\lambda$ based on the DRO formulation. A DRO problem aims to find a regression vector $\beta$ that minimizes the worst-case loss over an uncertainty set, namely
$\sup_{\mathbb{P} \in \mathcal{U}}\, {\rm E}_{\mathbb{P}} [ \ell(X,Y;\beta)]$.
Here, $\ell$ denotes a loss function, ${\rm E}_{\mathbb{P}}[\cdot]$ denotes the expectation with respect to a probability distribution $\mathbb{P}$ of $(X,Y)$, and $\mathcal{U}$ denotes the uncertainty set of probability measures.  The uncertainty set $\mathcal{U}$ specifies prior distributional information about $\mathbb{P}$, and it is usually constructed to include the unknown true  distribution with a probabilistic guarantee.
Let us now introduce some notation and give the DRO formulation related to \eqref{eq:square-root}. We denote the Dirac distribution at $(X_i,Y_i)$ by ${\bf 1}_{\{(X_i,Y_i)\}}$, the empirical distribution by $\mathbb{P}_N (dx,dy) := \frac{1}{N} \sum_{i=1}^{N} {\bf 1}_{\{(X_i,Y_i)\}}(dx,dy)$, and the squared-loss function by $ \ell (x,y;\beta) := (y - \beta^T x)^2$. Therefore, we have that ${\rm E}_{\mathbb{P}_N} [ \ell(X,Y;\beta)] = \frac{1}{N}\|Y - X\beta\|^2 $, and   \eqref{eq:square-root} takes an equivalent form
\begin{equation}\label{eq:square-root2}
	 \min_{\beta\in\mathbb{R}^n} \quad \left\{  \sqrt{ {\rm E}_{\mathbb{P}_N} [ \ell(X,Y;\beta)]} + \frac{\lambda}{\sqrt{N}} \cdot  p(\beta)\right\}.
\end{equation}
When $p(\beta)=\|\beta\|_q:= (\sum_{i=1}^n |\beta_i|^q)^{1/q}\,\,\forall\,\beta\in\mathbb{R}^n$ and $q\in [1,\infty)$, \eqref{eq:square-root2}  has the following equivalent DRO representation, as shown in \citep{blanchet2019robust},
\begin{equation}\label{eq:DRO}
	\inf_{\beta\in\mathbb{R}^n} \quad
	\sup_{\mathbb{P}}\Big\{ {\rm E}_{\mathbb{P}} [ \ell(X,Y;\beta)]
	\,\Big| \,
	\mathcal{D}_c (\mathbb{P},\,\mathbb{P}_N)\leq \frac{\lambda^2}{N}
	\Big\}
\end{equation}
where $\mathcal{D}_c (\mathbb{P},\mathbb{Q})$ is an optimal transport cost between two probability measures $\mathbb{P}$ and $\mathbb{Q}$ based on a suitably chosen cost function $c$ (see \eqref{eq:transport-cost} for the definition of $\mathcal{D}_c$), and   $\{\mathbb{P}\,|\,\mathcal{D}_c (\mathbb{P},\mathbb{P}_N)\leq \frac{\lambda^2}{N}\}$ is the uncertainty set centered at the empirical distribution $\mathbb{P}_N$  with radius $\frac{\lambda^2}{N}$. The inner maximization problem of \eqref{eq:DRO} accounts for all the probability measures that are  plausible variations of $\mathbb{P}_N$. The problem \eqref{eq:DRO} tries to minimize the worst-case loss and consequently it is likely to perform uniformly well around the empirical distribution.
Furthermore, we can interpret from  \eqref{eq:DRO}   that the regularization parameter $\lambda$  fully quantifies the radius of the uncertainty set $ \left\{\mathbb{P} \mid  \mathcal{D}_c(\mathbb{P},\mathbb{P}_N) \leq \frac{\lambda^2}{N}\right\}$.
In fact, when $p$ is the $\ell_q$ norm, \citet{blanchet2019robust} proved that an associated cost function $c$ defined by the dual norm of the $\ell_q$ norm will give rise to the equivalence between \eqref{eq:square-root2} and \eqref{eq:DRO}.
Besides, \citet{blanchet2017distributionally}  provided a DRO representation for the square-root group Lasso model, i.e., \eqref{eq:square-root} with $p$ being the group Lasso penalty function. The problem \eqref{eq:square-root} also admits a DRO representation when the regularizer is given by $p(\beta) = \|\beta\|_{\Lambda} := \sqrt{\beta^T \Lambda \beta}\,\,\forall\,\beta\in\mathbb{R}^n$ with $\Lambda$ being a given symmetric positive definite matrix, as shown in \citep[Theorem~1]{blanchet2019data}. Additionally, the connections between regularization and robust optimization have been extensively studied in the literature, e.g., \citep{el1997robust,xu2010robust,xu2009robustness,shafieezadeh2015distributionally,bertsimas2018characterization}.

One of our contributions in this paper is to provide a DRO representation for a general class of square-root regularized models  where the regularizer can take the general form as the sum of a simple norm and a seminorm. This broad class of regularizers can include most of the popular penalty functions in statistics and machine learning. As a key ingredient in the DRO formulation, the optimal transport cost
is given by a certain dual form of the regularizer. Despite the superior statistical properties, square-root regularized models are more complex and challenging to solve compared to their squared-loss counterparts because now the loss function and the penalty function
are both nonsmooth. Thus designing efficient algorithms that are capable of
solving high-dimensional square-root regularized models is an important task
for making these models practically useful. Here we develop a
proximal point dual semismooth Newton algorithmic framework for solving a generic
square-root regularized model. Specifically, we illustrate how the general framework
can be adopted to solve the square-root sparse group Lasso and  the square-root fused Lasso models, i.e.,
 for solving \eqref{eq:square-root} when $p$ is one of the following regularizers:
\begin{itemize}
\item the sparse group Lasso regularizer \citep{friedman2010note}
\begin{equation}\label{group-reg}
p(\beta) = w_1 \|\beta\|_1 + w_2\sum_{l=1}^g \omega_l \|\beta_{G_l}\|\,\,\forall\,\beta\in\mathbb{R}^n,
\end{equation}
where $w_1$ and $w_2$ are nonnegative regularization parameters,
$\omega_l > 0$ is the weight for the $l$th group, $G_l$'s form a partition of $\{1,\dots,n\}$, and $\beta_{G_l}$ is the  subvector of $\beta$ restricted
to $G_l$;

\item the fused Lasso regularizer \citep{tibshirani2005sparsity}
\begin{equation}\label{fused-reg}
p(\beta) = w_1\|\beta\|_1 + w_2\sum_{i=1}^{n-1}|\beta_i - \beta_{i+1}|\,\,\forall\,\beta\in\mathbb{R}^n,
\end{equation}
where $w_1$ and $w_2$ are nonnegative regularization parameters.
\end{itemize}

We should mention that although the square-root regularized models of the form \eqref{eq:square-root} enjoy nice statistical properties and insightful DRO interpretation, there is currently no efficient unified algorithmic framework for solving such models, possibly due to challenge posed by the nonsmoothness of the square-root loss function in addition to the nonsmoothness of the penalty function.
In fact, existing algorithms are limited to solving special cases such as the square-root Lasso model, and they are not efficient enough for solving large-scale problems.
\citet{belloni2011square} reformulated the square-root Lasso model into a second order cone programming (SOCP), and then applied TFOCS \citep{becker2011templates} (first order conic solvers), SDPT3 \citep{toh1999sdpt3,tutuncu2003solving} (interior point methods), and coordinatewise methods for solving the reformulated SOCP.
However, this SOCP reformulation will at least double the number of variables---thereby increasing the computational cost substantially, and those off-the-shelf methods are not efficient enough for solving large-scale problems. Additionally, an alternating direction method of multipliers  was applied in \citep{li2015flare} for solving the square-root Lasso model. But this approach may not be efficient for large-scale problems, as can be seen from the numerical experiments in \citep{tang2020sparse}.
Moreover, proximal gradient descent and proximal Newton methods were proposed in \citep{li2020fast} for solving the square-root Lasso model.
For solving the square-root group Lasso model,  \citet{bunea2013group} proposed a scaled thresholding-based iterative selection procedure. This method highly depends on the choice of the scaling parameter, and it usually requires many iterations for solving large-scale problems, as shown later in our numerical experiments. Moreover, we found that the convergence of the method  \citep[Theorem~3.1]{bunea2013group} is based on a regularity condition that depends on the sequence of iterates generated by the method and it may not hold generally.
In a recent work, \citet{tang2020sparse} proposed a semismooth Newton (SSN) based proximal majorization-minimization algorithm for solving nonconvex square-root-loss regression problems, and they demonstrated its efficiency for solving the square-root Lasso model.  To the best of our knowledge, an efficient and robust algorithm for solving the square-root regularized model that allows for a general regularizer is still in great demand but not yet available.
Inspired by the algorithm of \citep{tang2020sparse},  we propose in this paper a proximal point dual semismooth Newton algorithm (PPDNA) for solving \eqref{eq:square-root} with regularizer $p$ having the property that
its proximal mapping and the associated generalized Jacobian can be computed efficiently. In particular,
we implement the PPDNA algorithmic framework to solve
the square-root regularized model \eqref{eq:square-root} with the sparse group Lasso regularizer \eqref{group-reg} or the fused Lasso regularizer \eqref{fused-reg}.
In contrast to the algorithm in \citep{tang2020sparse} which does not impose structured sparsity, our algorithm will incorporate structured sparsity imposed by the sparse group Lasso or fused Lasso regularizer.
We should add that our PPDNA framework, just like the one in \citep{tang2020sparse}, is inspired
by the highly efficient SSN based proximal point algorithmic (PPA) framework developed
solving many squared-loss Lasso-type problems \citep{li2018efficiently,li2018highly,lin2019efficient,luo2019solving,zhang2020efficient}.
The key difference between our current framework and the one developed for a squared-loss
Lasso-type problem is that the latter is applied
to the dual problem having an essentially smooth strongly convex term in the objective function, but such a property is not present for our square-root regularized models.
Fortunately, as we shall see later, the desirable properties of the PPA and SSN, such as the fast convergence speed of the PPA and the ability to exploit second order sparsity in the SSN method, are preserved in our current PPDNA framework when there is no overfitting of the data.

The remaining parts of the paper are organized as follows.  In Section~\ref{sec:DRO} we show that the square-root regularized model \eqref{eq:square-root} has an equivalent DRO representation.
We design a  proximal point dual semismooth Newton algorithm
for solving  the square-root regularized model \eqref{eq:square-root} in Section~\ref{sec:algorithm}, and we elucidate the computational details of two cases where the regularizer is chosen to be the sparse group Lasso regularizer \eqref{group-reg} or the fused Lasso regularizer \eqref{fused-reg}. We conduct numerical experiments on synthetic and real data sets in  Section~\ref{sec:numerical-result}, and
give the conclusion in Section~\ref{sec:conclusion}.

\bigskip
\noindent{\bf Notation}
Throughout the paper any vector is understood to be a column vector.
We denote the inner product of two vectors $x$ and $y$ in $\mathbb{R}^n$ by $\langle x,y\rangle := x^T y$. We denote the Euclidean norm by $\|\cdot\|$ and the unit ball of Euclidean norm by $\mathcal{B}:= \left\{ x \in\mathbb{R}^n \mid \norm{x} \leq 1 \right\} $. For a matrix $X\in\mathbb{R}^{m\times n}$, we denote the operator norm of $X$ by $\norm{X}:=\inf\left\{c\mid \norm{Xv}\leq c\norm{v} \forall v\in\mathbb{R}^n \right\} $. We denote the elementwise multiplication by $\odot$.
We adopt the conventions of extended arithmetic, whereby $\infty\cdot 0 = 0 \cdot \infty = 0/0 = 0$ and $\infty - \infty = -\infty + \infty = 1/0 = \infty$.
For any subset $A\subseteq \mathbb{R}^n$, we denote the cardinality of $A$ as $|A|$.
We denote the vector of all ones by  $\bm{1}$ and the identity matrix by $I$. We denote the vector (or the matrix) of all zeros by $\bm{0}$.
For a matrix $X\in\mathbb{R}^{m\times n}$ and an index set $G\subseteq \{1,2,\dots,n\}$, we denote the range space (resp. null space) of $X$ by $\mathtt{Range}(X)$ (resp. $\mathtt{Null}(X)$) and the submatrix formed by the columns of $X$ corresponding to $G$ by $X_{G}$.
We denote the square diagonal matrix with the elements of vector $v$ on the main diagonal by $\mathtt{Diag}(v)$.
For a vector $\beta\in\mathbb{R}^n$, we denote the positive part of $\beta$ by $\beta^+:= \max\{\beta,0\} $, the vector with the signs of the corresponding elements of $\beta$ by $\mathtt{sign}(\beta)$ (the sign of a real number is $1$, $0$, or $-1$ if the number is positive, zero, or negative, respectively),
the restricted vector of $\beta$ to an index set $G\subseteq \{1,2,\dots,n\}$ by $\beta_{G}$, the support of $\beta$ by $\mathtt{supp}(\beta) := \{i\,|\,\beta_i \neq 0\}$.
Given a vector space $\mathcal{V}$, a norm on $\mathcal{V}$ is a nonnegative valued function $p:\,\mathcal{V} \to \mathbb{R}$ with the following properties:
for all $\lambda \in \mathbb{R}$ and $x,\,y \in \mathcal{V}$,
(1) (triangle inequality) $p(x + y) \leq p(x) + p(y)$,
(2) (absolute homogeneity) $p(\lambda x) = |\lambda| p(x)$,
(3) (positive definiteness) $p(x) = 0$ implies $x=0$.
A seminorm on $\mathcal{V}$ is a function $p:\,\mathcal{V} \to \mathbb{R}$ with the properties (1) and (2) above.
For a seminorm $p$ on $\mathbb{R}^n$, we define $p_*:\,\mathbb{R}^n \to [0,+\infty] $ by $p_*(y) := \sup_{x}\,\{\langle y,x \rangle \,|\, p(x) \leq 1\}\,\,\forall\,y\in\mathbb{R}^n$ (if $p$ is a norm then $p_*$ is also a norm, called its dual norm).
For a closed proper convex function $f:\,\mathbb{R}^n\to (-\infty,+\infty]$, we denote its effective domain  by $\mathtt{dom}(f) := \{x\,|\,f(x) < +\infty\}$, its Fenchel conjugate by $f^*(x) := \sup_{y\in\mathbb{R}^n} \{ \langle y,x \rangle - f(y)\}$,
the proximal mapping of $f$ at $x$ by
$\mathtt{prox}_{f}(x):= \arg\min_{y}\{ f(y) +\frac{1}{2}\|x-y\|^2 \}$,
and the Moreau envelope \citep{moreau1965proximite,yosida1980functional} of $f$ at $x$ by $\mathtt{M}_{f}(x):= \min_{y}\{ f(y) +\frac{1}{2}\|x-y\|^2 \}$. The Moreau envelope is continuously differentiable with the gradient $\nabla \mathtt{M}_{f}(x) = x - \mathtt{prox}_{f}(x)\,\,\forall \, x$.
The multivariate normal distribution with mean vector $\mu$ and covariance matrix $\Sigma$
is denoted by $\mathcal{N}(\mu,\Sigma)$.

\section{DRO Formulation of  Squared-loss  Linear Regression Problems}\label{sec:DRO}

In this section,  we will show that the square-root regularized model \eqref{eq:square-root} is equivalent
to a DRO formulation of a squared-loss linear regression problem, namely,
\begin{equation}\label{eq:DRO2}
	\inf_{\beta\in\mathbb{R}^n} \,\, \sup_{\mathbb{P}:\,\mathcal{D}_c (\mathbb{P},\mathbb{P}_N)\leq \delta}\, {\rm E}_{\mathbb{P}} [ \ell(X,Y;\beta)].
\end{equation}
In the above, $c$ is an appropriate optimal transport cost function
that will be defined explicitly later, and $\delta$ can be regarded as the radius of the uncertainty set centered at $\mathbb{P}_N$. Later we will show that the penalty parameter $\lambda$ in \eqref{eq:square-root} will fully quantify the radius $\delta$.

\subsection{Optimal Transport Costs}
We introduce in this section some notation and recall the optimal transport cost between probability measures; see \citep[Chapter~6]{villani2008optimal} for more details. For any two probability measures $\mathbb{P}$ and $\mathbb{Q}$ in $\mathbb{R}^{n+1}$, $\Pi(\mathbb{P},\mathbb{Q})$ denotes the set of all joint probability measures on $\mathbb{R}^{n+1} \times \mathbb{R}^{n+1}$ whose marginals are $\mathbb{P}$ and $\mathbb{Q}$.
For a given cost function $c:\,\mathbb{R}^{n+1} \times \mathbb{R}^{n+1} \to [0,\infty]$, where $c(u,v)$ is the cost for transporting one unit of mass from $u$ to $v$,
 the optimal transport cost between $\mathbb{P}$ and $\mathbb{Q}$ is defined as
\begin{equation}\label{eq:transport-cost}
  	\mathcal{D}_{c}(\mathbb{P},\mathbb{Q}) := \inf_{\pi \in \Pi(\mathbb{P},\mathbb{Q})}  \int c(u,v) d\pi(u,v).
\end{equation}
We assume that $c(u,u)=0$ for any $u \in \mathbb{R}^{n+1}$. For any nonnegative lower semicontinuous cost function $c$, it is shown in \citep[Theorem~4.1]{villani2008optimal} that the infimum in \eqref{eq:transport-cost} is attainable. Intuitively, one can regard \eqref{eq:transport-cost} as a kind of distance between two measures $\mathbb{P}$ and $\mathbb{Q}$, but strictly speaking, it is not guaranteed to satisfy the axioms of a distance. However, one can obtain a distance from \eqref{eq:transport-cost} when the cost is defined in terms of a distance. For example, if the cost function is defined by the $\ell_q$-norm as
$c(u,v) = \|u - v\|_q^{\rho},\,q \geq 1,\, \rho \geq 1$, then $\mathcal{D}_c(\cdot,\cdot)$ is the well known Wasserstein distance of order $\rho$ (also known as the optimal transport distance or the earth mover's distance).
Wasserstein distances and Wasserstein barycenters have recently become
very popular and are widely applied in many applications \citep{yang2021fast,bigot2018characterization,cuturi2014fast,li2008real,ye2017fast,rabin2011wasserstein}.   In contrast to standard Wasserstein distances, the cost function $c$ in our analysis are more general in that we allow for  lower semicontinuous  cost functions that may take infinite values. As we will see in Theorem~\ref{thm:main1}, a judicious choice of $c$ will give rise to the equivalence between \eqref{eq:square-root} and \eqref{eq:DRO2}.

\subsection{DRO Formulation}
In this section, we will give a DRO representation of the square-root regularized model \eqref{eq:square-root} for a broad class of regularizers
$p$
expressed as the sum of a simple norm and a seminorm as follows.
Let $P:\,\mathbb{R}^n \to \mathbb{R}$ and $Q:\,\mathbb{R}^s \to \mathbb{R}$ be two given norms. We consider $p:\,\mathbb{R}^n \to \mathbb{R}$ defined by
\begin{equation}\label{eq:pform}
p(\beta):= w_1 P(\beta) + w_2 Q(B \beta),\,\beta\in\mathbb{R}^n,
\end{equation}
where $B \in \mathbb{R}^{s\times n}$ is a given matrix, and $w_1$ and $w_2$ are nonnegative scalars adding up to one.
Such a regularizer $p$ and the corresponding dual function $p_*$ are essential ingredients in the definition of the optimal transport cost in the DRO formulation. In  Proposition~\ref{prop:pnorm}, we give a relatively explicit form of $p_*$.

\begin{proposition}\label{prop:pnorm}
Let $p:\mathbb{R}^n \to \mathbb{R}$ be the function defined
in \eqref{eq:pform}. Then the function
\begin{equation}\label{eq:def-pstar}
p_*(\alpha) := \sup_{\beta}\,\{\langle \alpha,\beta \rangle \,|\, p(\beta) \leq 1\},\,\alpha\in\mathbb{R}^n
\end{equation}
admits the form
\begin{equation}\label{eq:pstar}
p_*(\alpha) =  \inf_{\tilde{\alpha}\in\mathbb{R}^n,\,\bar{\alpha}\in\mathbb{R}^s } {\sup_{t \in [0,1]}}
\left\{ t \frac{P_*(\tilde{\alpha})}{w_1} + (1-t) \frac{Q_*(\bar{\alpha})}{w_2}\,\bigg|\,
\tilde{\alpha} + B^T \bar{\alpha} = \alpha \right\},\,\alpha\in\mathbb{R}^n.
\end{equation}
\end{proposition}

Before presenting the proof, we would like to give an explanation of this proposition. By virtue of the convention of extended arithmetic ($0/0 = 0,\,1/0 = \infty$), the formula \eqref{eq:pstar} can also include the special cases with $w_1=0$ or $w_2=0$, namely,
\begin{equation}\label{eq:special}
p_*(\alpha) = \left\{
\begin{array}{ll}
\displaystyle P_*(\alpha), & \mbox{if } w_1 = 1,\,w_2 = 0,\\[6pt]
\displaystyle \inf_{\bar{\alpha}\in\mathbb{R}^s} \{Q_*(\bar{\alpha})\,|\,B^T \bar{\alpha} = \alpha\}, & \mbox{if } w_1 = 0,\,w_2 = 1.
\end{array}
\right.
\end{equation}
The convention of extended arithmetic allows  us  to have a uniform expression \eqref{eq:pstar} without having to separately write out different cases. Additionally, one can see that $p_*$ in the second case of \eqref{eq:special} will have finite values if and only if $\alpha \in \mathtt{Range}(B^T)$. Note that in this case  $p(\cdot) = Q(B \cdot)$ may only be a seminorm on $\mathbb{R}^n$ and therefore $p_*$ may take infinite values. We will further characterize the properties of $p_*$ when $p(\cdot) = Q(B \cdot)$ is a seminorm later in Proposition~\ref{prop:seminorm}.

\begin{proof}
For positive coefficients $w_1 > 0$ and $w_2 > 0$,
take an arbitrary $\bm{0} \neq \alpha \in \mathbb{R}^n$. We have that
\begin{equation*}
\begin{array}{l}
\displaystyle p_*(\alpha) = \sup_{\beta}\,\{\langle \alpha,\beta \rangle \,|\, p(\beta) \leq 1\}  \\[6pt]
\displaystyle= \sup_{\beta} \inf_{u\geq 0} \,\{\langle \alpha,\beta \rangle + u(1 - p(\beta))\}  \\[6pt]
\displaystyle = \inf_{u\geq 0} \sup_{\beta}\,\{\langle \alpha,\beta \rangle + u - uw_1P(\beta) - uw_2Q(B\beta)\}  \\[6pt]
\displaystyle= \inf_{u\geq 0} \sup_{\beta}\,\{\langle \alpha,\beta \rangle + u - uw_1 \sup_{\tilde{\alpha}:\,P_*(\tilde{\alpha})\leq1} \langle \beta,\tilde{\alpha}\rangle - uw_2 \sup_{\bar{\alpha}:\,Q_*(\bar{\alpha}) \leq 1} \, \langle B\beta,\bar{\alpha} \rangle \}  \\[6pt]
\displaystyle = \inf_{u \geq 0} \, \{ u + \sup_{\beta} \inf_{\tilde{\alpha}:\,P_*(\tilde{\alpha})\leq1,\,\bar{\alpha}:\,Q_*(\bar{\alpha}) \leq 1} \, \langle \beta, \alpha - u(w_1\tilde{\alpha} +w_2 B^T \bar{\alpha})\rangle \}  \\[6pt]
\displaystyle = \inf_{u \geq 0} \, \{ u + \inf_{\tilde{\alpha}:\,P_*(\tilde{\alpha})\leq1,\,\bar{\alpha}:\,Q_*(\bar{\alpha}) \leq 1}\sup_{\beta}  \, \langle \beta, \alpha - u(w_1\tilde{\alpha} +w_2 B^T \bar{\alpha})\rangle \} .
\end{array}
\end{equation*}
The third equality follows from  the strong duality theorem \citep[Theorem~28.2 and 28.4]{rockafellar1997convex}; the fourth equality applies the definition of a dual norm; and the last equality follows from a standard minimax theorem \citep[Corallary~3.3]{sion1958general}. The supremum over $\beta$ renders that $\alpha - u(w_1\tilde{\alpha} +w_2 B^T \bar{\alpha}) =\bm{0}$. Therefore, we have that
\begin{equation*}
\begin{array}{l}
\displaystyle \inf_{u \geq 0} \, \{ u + \inf_{\tilde{\alpha}:\,P_*(\tilde{\alpha})\leq1,\,\bar{\alpha}:\,Q_*(\bar{\alpha}) \leq 1}\sup_{\beta}  \, \langle \beta, \alpha - u(w_1\tilde{\alpha} +w_2 B^T \bar{\alpha})\rangle \} \\[6pt]
\displaystyle  =	{\inf_{u,\,\tilde{\alpha},\,\bar{\alpha}}}\, \{u\,|\, u\geq 0,\,P_*(\tilde{\alpha})\leq 1,\,Q_*(\bar{\alpha}) \leq 1,\,u(w_1\tilde{\alpha} +w_2 B^T \bar{\alpha}) = \alpha\}  \\[6pt]
\displaystyle = {\inf_{u,\,\tilde{\alpha},\,\bar{\alpha}}}\, \{u\,|\, u > 0,\,P_*(\tilde{\alpha})\leq uw_1,\,Q_*(\bar{\alpha}) \leq uw_2,\, \tilde{\alpha} + B^T \bar{\alpha} = \alpha\}  \\[6pt]
\displaystyle = {\inf_{\tilde{\alpha},\,\bar{\alpha}}}\, \{\max ( P_*(\tilde{\alpha})/w_1,  Q_*(\bar{\alpha})/w_2 ) \,|\,   \tilde{\alpha} + B^T \bar{\alpha} = \alpha\}  \\[6pt]
\displaystyle = \inf_{\tilde{\alpha},\,\bar{\alpha}}{\sup_{t \in [0,1]}}
\left\{ t \frac{P_*(\tilde{\alpha})}{w_1} + (1-t) \frac{Q_*(\bar{\alpha})}{w_2}\,\bigg|\,
\tilde{\alpha} + B^T \bar{\alpha} = \alpha \right\},
\end{array}
\end{equation*}
where the second equality holds since $\alpha \neq \bm{0}$, $w_1>0$, $w_2>0$, and we can simply replace the variable $\tilde{\alpha}$ by $\tilde{\alpha}/(uw_1)$ and the variable $\bar{\alpha}$ by $\bar{\alpha}/(uw_2)$. Therefore, \eqref{eq:pstar} holds for $w_1>0$ and $w_2>0$.

For the case when $w_1 = 1$ and $w_2 = 0$, the required result
follows trivially.

Next, we prove the result for case when $w_1 = 0$ and $w_2 = 1$, where $p(\cdot) = Q(B\, \cdot)$ is a seminorm but it is not necessarily  a norm. Take an arbitrary $\bm{0} \neq \alpha \in \mathbb{R}^n$. We have that
\begin{equation*}
\begin{array}{l}
\displaystyle p_*(\alpha) = \sup_{\beta}\,\{\langle \alpha,\beta \rangle \,|\, Q(B\beta) \leq 1\}
\displaystyle= \sup_{\beta} \inf_{u\geq 0} \,\{\langle \alpha,\beta \rangle + u(1 - Q(B\beta))\}  \\[6pt]
  \displaystyle \leq \inf_{u\geq 0} \sup_{\beta}\,\{\langle \alpha,\beta \rangle + u(1 - Q(B\beta))\} 
  \\[6pt]
\displaystyle= \inf_{u\geq 0} \sup_{\beta}\,\{\langle \alpha,\beta \rangle + u - u \sup_{\bar{\alpha}:\,Q_*(\bar{\alpha}) \leq 1} \, \langle B\beta,\bar{\alpha} \rangle \}  \\[6pt]
\displaystyle = \inf_{u \geq 0} \, \{ u + \sup_{\beta} \inf_{\bar{\alpha}:\,Q_*(\bar{\alpha}) \leq 1} \, \langle \beta, \alpha - u B^T \bar{\alpha}\rangle \}  \\[6pt]
\displaystyle = \inf_{u \geq 0} \, \{ u + \inf_{\bar{\alpha}:\,Q_*(\bar{\alpha}) \leq 1} \sup_{\beta} \, \langle \beta, \alpha - u B^T \bar{\alpha}\rangle \}.
\end{array}
\end{equation*}
The third inequality follows from  the exchange of sup and inf \citep[Lemma~36.1]{rockafellar1997convex}; the fourth equality applies the definition of a dual norm; and the last equality follows from a standard minimax theorem \citep[Corallary~3.3]{sion1958general}. The supremum over $\beta$ renders that $\alpha - u B^T \bar{\alpha} = {\bm{0}}$. Therefore, we have that
\begin{equation*}
\begin{array}{l}
\displaystyle  \inf_{u \geq 0} \, \{ u + \inf_{\bar{\alpha}:\,Q_*(\bar{\alpha}) \leq 1} \sup_{\beta} \, \langle \beta, \alpha - u B^T \bar{\alpha}\rangle \}
\displaystyle  =	\inf_{u,\,\bar{\alpha}}\, \{u\,|\, u\geq 0,\,Q_*(\bar{\alpha}) \leq 1,\,u B^T \bar{\alpha} = \alpha\}  \\[6pt]
\displaystyle = \inf_{u,\,\bar{\alpha}}\, \{u\,|\, u > 0,\,Q_*(\bar{\alpha}) \leq u,\, B^T \bar{\alpha} = \alpha\}
\displaystyle = \inf_{\bar{\alpha}}\, \{Q_*(\bar{\alpha})\,|\,B^T \bar{\alpha} = \alpha \},
\end{array}
\end{equation*}
where the second equality holds since $\alpha \neq {\bm{0}}$ and we can simply replace the variable $\bar{\alpha}$ by $\bar{\alpha}/u$.

When ${\bm{0}}\neq \alpha \in \mathtt{Range}(B^T)$,  $\inf_{\bar{\alpha}}\, \{Q_*(\bar{\alpha})\,|\,B^T \bar{\alpha} = \alpha \}$ is finite. This, together with the strong duality theorem \citep[Theorem~28.2 and 28.4]{rockafellar1997convex}, implies that equality holds in the above inequality ($\ast$). Namely, $p_*(\alpha) = \inf_{\bar{\alpha}}\, \{Q_*(\bar{\alpha})\,|\,B^T \bar{\alpha} = \alpha \} < +\infty$ for  ${\bm{0}}\neq \alpha \in \mathtt{Range}(B^T)$.
	
When ${\bm{0}}\neq \alpha \notin \mathtt{Range}(B^T)$,  $\inf_{\bar{\alpha}}\, \{Q_*(\bar{\alpha})\,|\,B^T \bar{\alpha} = \alpha \}$ is infinite.	We can have an orthogonal decomposition of $\alpha$: $\alpha = \alpha_r + \alpha_n,\,\alpha_r\in\mathtt{Range}(B^T),\,{\bm{0}}\neq \alpha_n \in \mathtt{Null}(B)$. Note that $p(\alpha_n)=0$, and by \eqref{eq:def-pstar} it holds that $p_*(\alpha) \geq \langle \alpha, k \alpha_n\rangle = k \|\alpha_n\|^2 \to +\infty$, as $k\to +\infty$. Therefore, $p_*(\alpha) = \inf_{\bar{\alpha}}\, \{Q_*(\bar{\alpha})\,|\,B^T \bar{\alpha} = \alpha \} = +\infty $ for  ${\bm{0}}\neq \alpha \notin \mathtt{Range}(B^T)$.

Together with $p_*(0)=0$, the proof is completed.
\end{proof}

For a seminorm $p$, we specify some properties of  $p$ and $p_*$. Note that if $p$ is a norm on $\mathbb{R}^n$, then $p_*$ is also a norm on $\mathbb{R}^n$ and the following properties (b) and (c) are standard.

\begin{proposition} \label{prop:seminorm}
Let $B \in \mathbb{R}^{s\times n}$ be a given matrix, $Q:\,\mathbb{R}^s \to \mathbb{R}$ be a norm,  and $p:\,\mathbb{R}^n \to \mathbb{R}$ be defined by
\begin{equation}\label{eq:pform2}
p(\beta) =  Q(B \beta),\,\beta\in\mathbb{R}^n.
\end{equation}
Then $p_*$ defined by \eqref{eq:def-pstar}  admits the expression
\begin{equation}\label{eq:pstar2}
p_*(\alpha) = \inf_{\bar{\alpha}\in\mathbb{R}^s} \{Q_*(\bar{\alpha})\,|\,B^T \bar{\alpha} = \alpha\},\,\alpha\in\mathbb{R}^n,
\end{equation}
and the following  holds:
\begin{itemize}
\item[{\rm(a)}] $\mathtt{dom}(p_*) = \mathtt{Range}(B^T)$ and $p_*:\,\mathtt{Range}(B^T)\to \mathbb{R}$ is a norm on the vector space $\mathtt{Range}(B^T)$;

\item[{\rm(b)}]   $\alpha^T\beta \leq  p(\beta) p_*(\alpha)$  $\forall\,\alpha,\beta\in\mathbb{R}^n$;

\item[{\rm(c)}] $p(\beta) = \sup_{\alpha} \{\alpha^T\beta\,|\,p_*(\alpha) \leq 1\}$ $\forall \, \beta\in\mathbb{R}^n$. Moreover, this supremum is achievable. Namely, there exists ${\alpha}_{\beta}\in\mathtt{Range}(B^T) $ such that $p_*({\alpha}_{\beta})=1$ and $p(\beta) = {\alpha}_{\beta}^T\beta $.
\end{itemize}
\end{proposition}

\begin{proof}
\eqref{eq:pstar2} follows from Proposition~\ref{prop:pnorm}, and it implies that $\mathtt{dom}(p_*) = \mathtt{Range}(B^T)$.
First, we prove that $p_*:\,\mathtt{Range}(B^T)\to \mathbb{R}$ is a norm on $\mathtt{Range}(B^T)$. The triangle inequality is inherited from that of $Q_*$:
$ p_*(\alpha + \beta) = \inf_{\gamma} \{Q_*(\gamma)\,|\, B^T \gamma = \alpha + \beta\}
\leq \inf_{\bar{\alpha}, \bar{\beta}} \{Q_*(\bar{\alpha} + \bar{\beta})\,|\,B^T\bar{\alpha} = \alpha,B^T\bar{\beta}=\beta\}
\leq \inf_{\bar{\alpha}} \{Q_*(\bar{\alpha})\,|\,B^T\bar{\alpha} = \alpha\}
+  \inf_{\bar{\beta}} \{Q_*(\bar{\beta})\,|\,B^T\bar{\beta}=\beta\}
= p_*(\alpha) +  p_*(\beta) $.  The absolute homogeneity is obvious from that of  $Q_*$. If $p_*(\alpha) = 0$, then by the property of infimum there exist a sequence $\bar{\alpha}_k,k=1,2,\dots$ such that $B^T\bar{\alpha}_k = \alpha$ and $Q_*(\bar{\alpha}_k) < 1/k$. Since $Q_*$ is a norm on $\mathbb{R}^s$, the latter
further implies that $\bar{\alpha}_k\to 0$ and hence $\alpha=0$. Therefore, $p_*$ is a norm on $\mathtt{Range}(B^T)$.

Second, we prove (b). For any $\alpha$, $\beta$, and $\bar{\alpha}$ satisfying $B^T\bar{\alpha} = \alpha$, we have the inequality $\alpha^T \beta = (B\beta)^T\bar{\alpha} \leq Q(B\beta) Q_*(\bar{\alpha})$ since $Q$ is a norm. We can deduce (b) by taking infimum  over $\bar{\alpha}$.

Third, we prove (c). On one hand, it follows from (b) that $\sup_{\alpha} \{\alpha^T\beta\,|\,p_*(\alpha) \leq 1\} \leq \sup_{\alpha} \{p(\beta)p_*(\alpha)\,|\,p_*(\alpha) \leq 1\} \leq p(\beta)\,\,\forall\,\beta$.
On the other hand, we have that
$\sup_{\alpha} \{\alpha^T\beta\,|\,p_*(\alpha) \leq 1\}
= \sup_{\alpha\in \mathtt{Range}(B^T)} \{\alpha^T\beta\,|\, \inf_{\bar{\alpha}}\{Q_*(\bar{\alpha})\,|\,B^T\bar{\alpha} = \alpha\}\leq 1\}
= \sup_{\gamma} \{(B\beta)^T\gamma\,|\, \inf_{\bar{\alpha}}\{Q_*(\bar{\alpha})\,|\,B^T\bar{\alpha} = B^T \gamma\}\leq 1\}
\geq \sup_{\gamma} \{(B\beta)^T\gamma\,|\, Q_*(\gamma) \leq 1\}
= Q(B\beta)=p(\beta)
$,
where the second last equality follows from the properties of the norm $Q$ and its dual norm $Q_*$. Moreover, from (a) we know that $p_*$ is a norm on $\mathtt{Range}(B^T)$ and therefore the set $\{\alpha\,|\,p_*(\alpha) \leq 1\}$ is compact. Then it is easy to show that  this supremum is achievable  at the boundary. The proof is completed.
\end{proof}

Proposition~\ref{prop:seminorm} includes a similar result in \citep{maurer2012structured} about the operator norm $\|\beta\|_{\mathcal{M}^*} := \sup_{1\leq l \leq g} \left\{ \|M_{(l)}\beta\|  \right\}\,\,\forall \,\beta \in \mathbb{R}^n$ and its dual form, with $\mathcal{M}:= \{M_{(l)}\}_{1\leq l \leq g}$ being a set of symmetric matrices $M_{(l)} \in\mathbb{S}^n$. See Appendix~\ref{sec:remark} for details.

The next proposition about strong duality  is a direct application of \citep[Theorem~1]{blanchet2019quantifying}. It shows that the inner maximization in the DRO problem \eqref{eq:DRO2} has a nice univariate dual problem.
\begin{proposition}\citep[Proposition~1]{blanchet2019robust}\label{prop:blanchet}
Let $c: \, \mathbb{R}^{n+1} \times \mathbb{R}^{n+1} \to [0, +\infty] $ be a lower semicontinuous cost function satisfying $c ((x, y),(x', y'))=0 $ whenever $(x, y)=(x',y')$. For $ \gamma \geq 0 $ and a loss function $ \ell(x, y ; \beta)$ that is upper semicontinuous in $ (x, y) $ for each $ \beta $, define
\begin{equation*}
\phi_{\gamma} (X_{i}, Y_{i} ; \beta ):=\sup_{u \in \mathbb{R}^{n},\, v \in \mathbb{R}}\,
\left\{ \ell(u, v ; \beta)-\gamma c((u, v),(X_{i}, Y_{i})) \right\}.
\end{equation*}
Then
\begin{equation*}
\sup_{\mathbb{P}:\, \mathcal{D}_{c} (\mathbb{P}, \mathbb{P}_N ) \leq \delta}\, {\rm E}_{\mathbb{P}}[\ell(X, Y ; \beta)]
= \inf_{\gamma \geq 0} \, \left\{\gamma \delta + \frac{1}{N} \sum_{i=1}^{N} \phi_{\gamma} (X_{i}, Y_{i} ; \beta ) \right \}.
\end{equation*}
Consequently, the DRO problem \eqref{eq:DRO2} reduces to
\begin{equation*}
\inf_{\beta\in\mathbb{R}^n}\,\, \sup_{\mathbb{P}:\, \mathcal{D}_c(\mathbb{P},\mathbb{P}_N) \leq \delta } \, {\rm E}_{\mathbb{P}}  [\ell (X,Y;\beta)  ]
=\inf _{\beta \in \mathbb{R}^{n}}\,\, \inf _{\gamma \geq 0}\, \left\{\gamma \delta+\frac{1}{N} \sum_{i=1}^{N} \phi_{\gamma} (X_{i}, Y_{i} ; \beta ) \right\}.
\end{equation*}
\end{proposition}

Based on Proposition~\ref{prop:blanchet}, we will prove in the next main theorem the equivalence between the square-root regularized model \eqref{eq:square-root} and the DRO formulation \eqref{eq:DRO2}  by finding an explicit form of $\phi_{\gamma}$. This theorem not only unifies existing results \citep{blanchet2017distributionally,blanchet2019robust} but also include a broader class of regularizers; it is applicable for any regularizer in the
additive form expressed in \eqref{eq:pform}.
\begin{theorem} \label{thm:main1}
Consider the {squared-loss} function $ \ell (x,y;\beta) = (y - \beta^T x)^2$ and  a regularizer $p$ of the form \eqref{eq:pform}. Let the cost function $c: \, \mathbb{R}^{n+1} \times \mathbb{R}^{n+1} \to [0, +\infty] $ be defined by
\begin{equation*}
c ((u,v),(x,y)):= \left\{
\begin{array}{ll}
\left(p_*(u-x)\right)^2, & \mbox{if } v = y, \\[6pt]
+\infty, & \mbox{otherwise},
\end{array} \right.
\end{equation*}
and the associated optimal transport cost $\mathcal{D}_c(\cdot,\cdot)$ be defined by \eqref{eq:transport-cost}.	
Then it holds that
\begin{equation*}
\inf_{\beta\in\mathbb{R}^n} \,\,\sup_{\mathbb{P}:\, \mathcal{D}_c(\mathbb{P},\mathbb{P}_N) \leq \delta }  \,{\rm E}_{\mathbb{P}} [\ell (X,Y;\beta) ]
=\frac{1}{N} \cdot \inf_{\beta \in \mathbb{R}^{n}} \, \left\{ \|Y - X\beta\| + \sqrt{\delta N} \, p(\beta) \right\}^2.
\end{equation*}
\end{theorem}

\begin{proof}
Since $c$ is lower semicontinuous with
$c((u,v),(u,v))=0$ for any $(u,v)$,  and  $\ell(x,y;\beta)$ is upper semicontinuous in $(x,y)$ for each $\beta$,
based on Proposition~\ref{prop:blanchet}, we can prove the  required result via finding an explicit form  of $\phi_{\gamma}$. Take $\gamma \geq 0$ and $\beta\in\mathbb{R}^n$ arbitrarily.
By the definitions of $\phi_{\gamma},\,c$, and $\ell$, we have that
\begin{equation*}
\phi_\gamma(X_i,Y_i;\beta) = \sup_{u}\, \left\{ (Y_i - \beta^Tu)^2 - \gamma (p_*(u-X_i))^2 \,|\, u-X_i \in \mathtt{dom}(p_*)\right\}.
\end{equation*}
For notational simplicity, we denote  $ \Delta:= u-X_i $ and $ Z_i = Y_i - \beta^TX_i $. Then it holds that
\begin{align}
\phi_\gamma(X_i,Y_i;\beta) & =
\sup_{\Delta\in \mathtt{dom}(p_*)}\, \left\{ Z_i^2 - 2 Z_i \beta^T\Delta  + (\beta^T\Delta)^2  - \gamma (p_*(\Delta))^2 \right\} \label{pf:2} \\
&  \leq Z_i^2 + \sup_{\Delta\in \mathtt{dom}(p_*)}\, \left\{  2 |Z_i| p(\beta)p_*(\Delta) + \left((p(\beta))^2-\gamma\right)(p_*(\Delta))^2  \right\}. \label{pf:3}
\end{align}
By the property of a norm and Proposition~\ref{prop:seminorm}, there exists $\alpha_{\beta} \in \mathtt{dom}(p_*)$ such that $p_*(\alpha_{\beta}) = 1$ and $p(\beta) = \alpha_{\beta}^T \beta$.

If $(p(\beta))^2 < \gamma$, then the quadratic function in terms of $p_*(\Delta)$ in the above supremum problem \eqref{pf:3} is bounded  by the value $\frac{Z_i^2 (p(\beta))^2}{\gamma - (p(\beta))^2}$ at the stationary point $ p_*(\Delta) = \frac{|Z_i| p(\beta)}{\gamma - (p(\beta))^2}$. That is,
$ \phi_\gamma(X_i,Y_i;\beta) \leq Z_i^2 + \frac{Z_i^2 (p(\beta))^2}{\gamma - (p(\beta))^2} $.
It is easy to check that this equality  is achievable when one substitutes $\Delta = \frac{-Z_i p(\beta)}{\gamma - (p(\beta))^2} \alpha_{\beta}$ into \eqref{pf:2}.

If $(p(\beta))^2 > \gamma$, by taking $\Delta = k \alpha_{\beta}$, $k > 0$ in \eqref{pf:2}, then it holds that
$
\phi_\gamma(X_i,Y_i;\beta) \geq Z_i^2 - 2 k Z_i p(\beta) + k^2 (p(\beta) - \gamma)
$.
Since $k$ can be arbitrarily large, we deduce that $\phi_\gamma(X_i,Y_i;\beta)=+\infty$.

Lastly, we consider the case where $(p(\beta))^2 = \gamma$. If $Z_i p(\beta) = 0$, we have that $\phi_\gamma(X_i,Y_i;\beta)= Z_i^2$ since \eqref{pf:2} implies that $\phi_\gamma(X_i,Y_i;\beta) \geq Z_i^2$ by taking {$\Delta=\bm{0}$} and \eqref{pf:3} implies that $\phi_\gamma(X_i,Y_i;\beta) \leq Z_i^2$. If $Z_i p(\beta) \neq 0$, by taking $\Delta = -k Z_i p(\beta) \alpha_{\beta}$ with $k > 0$ in \eqref{pf:2}, then it holds that
$
\phi_\gamma(X_i,Y_i;\beta) \geq Z_i^2 + 2 k Z_i^2 (p(\beta))^2
$.
Since $k$ can be arbitrarily large, we deduce that $\phi_\gamma(X_i,Y_i;\beta)=+\infty$.
Therefore, when $(p(\beta))^2 = \gamma$, we have that
$
\phi_\gamma(X_i,Y_i;\beta) = \left\{\begin{array}{ll}
Z_i^2, & \mbox{if } Z_i p(\beta)=0, \\
+\infty, & \mbox{if } Z_i p(\beta)\neq 0.
\end{array}\right.
$

In summary, we can have a unified expression if we adopt the conventions of extended arithmetic given in the paragraph on notation, i.e.,
\begin{equation}\label{pf:5}
\phi_\gamma(X_i,Y_i;\beta) = \left\{\begin{array}{ll}
\displaystyle Z_i^2 + \frac{Z_i^2 (p(\beta))^2}{\gamma - (p(\beta))^2}, & \mbox{if $(p(\beta))^2 < \gamma$, or
$(p(\beta))^2 = \gamma$ and $Z_i p(\beta)=0$},
\\[6pt]
+\infty, & \mbox{if  $(p(\beta))^2 > \gamma$, or $(p(\beta))^2 = \gamma$ and $Z_i p(\beta)\not=0$}.
\end{array}\right.
\end{equation}

Next, we consider the problem
\begin{equation}\label{pf:6}
\inf_{\gamma \geq 0} \, \left\{\gamma \delta + \frac{1}{N} \sum_{i=1}^{N} \phi_{\gamma} (X_{i}, Y_{i} ; \beta ) \right \}.
\end{equation}
If $\|Y - X\beta\| = 0$ (namely, $Z_i = 0,\,\forall\, i$), or $p(\beta) = 0$, then it follows from \eqref{pf:5} that \eqref{pf:6} reduces to
\begin{equation*}
\inf_{\gamma \geq (p(\beta))^2} \, \left\{\gamma \delta + \frac{1}{N} \sum_{i=1}^{N} Z_i^2 \right \} = \frac{1}{N} \|Y - X\beta\|^2 + \delta (p(\beta))^2 = \left( \frac{1}{\sqrt{N}}\|Y - X\beta\| + \sqrt{\delta} p(\beta) \right)^2.
\end{equation*}
If $\|Y - X\beta\| \neq 0$ and $p(\beta) \neq 0$, then it follows from \eqref{pf:5} that \eqref{pf:6} reduces to
\begin{equation*}
\inf_{\gamma > (p(\beta))^2} \, \left\{\gamma \delta + \frac{1}{N} \sum_{i=1}^{N} \left( Z_i^2 + \frac{Z_i^2 (p(\beta))^2}{\gamma - (p(\beta))^2}\right) \right \} = \left( \frac{1}{\sqrt{N}}\|Y - X\beta\| + \sqrt{\delta} p(\beta) \right)^2,
\end{equation*}
where the minimal value is achieved at  $\gamma = \frac{1}{\sqrt{N \delta}} \|Y - X\beta\| p(\beta) + (p(\beta))^2$. Thus far, we have that
$ \inf_{\gamma \geq 0} \, \left\{\gamma \delta + \frac{1}{N} \sum_{i=1}^{N} \phi_{\gamma} (X_{i}, Y_{i} ; \beta ) \right \}
=\left( \frac{1}{\sqrt{N}}\|Y - X\beta\| + \sqrt{\delta} p(\beta) \right)^2
$.
We further take minimization over $\beta$ on both sides. This, together with Proposition~\ref{prop:blanchet}, completes the proof.
\end{proof}

Theorem~\ref{thm:main1} is applicable for any regularizer of the additive form \eqref{eq:pform}. This form may be a seminorm when $p(\cdot) = Q(B\cdot)$. In this case, one may take into consideration the effective domain of $p_*$, and then the cost function can be written as
\begin{equation*}
c ((u,v),(x,y)):= \left\{
\begin{array}{ll}
\left(p_*(u-x)\right)^2, & \mbox{if } v = y \mbox{ and } u-x\in\mathtt{Range}(B^T), \\[6pt]
+\infty, & \mbox{otherwise}.
\end{array} \right.
\end{equation*}
Our equivalence result in Theorem~\ref{thm:main1} can cover a broad class of regularizers composed  of a norm and a seminorm as in \eqref{eq:pform}.
We can obtain from Theorem~\ref{thm:main1} that the following square-root regularized estimators have equivalent DRO formulations:
\begin{itemize}
\item the square-root Lasso estimator \citep{belloni2011square}; a solution of \eqref{eq:square-root} with $p(\beta) = \|\beta\|_1,\,\beta\in\mathbb{R}^n$;

\item the square-root sparse group Lasso estimator \citep[Section~4.4]{stucky2017sharp}; a solution of \eqref{eq:square-root} with $p$ given by \eqref{group-reg};

\item the square-root SLOPE estimator \citep[Section~4.3]{stucky2017sharp}; a solution of \eqref{eq:square-root} with $p$ being the following weighted and sorted $\ell_1$ norm with respect to a nonincreasing sequence of weights $\omega_1\geq \dots \geq \omega_n > 0$:
$
p(\beta) = \sum_{i=1}^n \omega_i |\beta_i^{\downarrow}|,
$
where $\beta^{\downarrow}$ is a vector obtained from $\beta$  by sorting its entries in nonincreasing order of magnitude. This $p$ was shown to be a norm in \citep[Lemma~2]{zeng2014decreasing};

\item the square-root fused Lasso estimator \citep{MR4159660}; a solution of \eqref{eq:square-root} with $p$ being given by \eqref{fused-reg};

\item the SSASR estimator \citep{xie2020structured}; a solution of \eqref{eq:square-root} with $p$ being given by
$
p(\beta) = w_1 \sum_{l=1}^g \sqrt{|G_l|} \|\beta_{G_l}\|
+ w_2 \sum_{l=1}^g \sqrt{|G_l|} \| B_l \beta_{G_l} \|,
$
where the matrix $B_l \in \mathbb{R}^{(|G_l| - 1) \times |G_l|} $ is defined as $B_l x = (x_1 - x_2, \dots, x_{|G_l| - 1} - x_{|G_l|})^T\,\, \forall \, x \in \mathbb{R}^{|G_l|}$ if $|G_l| \geq 2$; and $B_l = 0$ if $|G_l| = 1$.

\end{itemize}

\section{A Proximal Point Dual Semismooth Newton Algorithm for Solving the Square-root Regularized Problem} \label{sec:algorithm}

In this section we aim to design a fast algorithm to solve the square-root regularized problem \eqref{eq:square-root}.
It can be rewritten as follows with an auxiliary vector $y\in\mathbb{R}^N$:
\begin{equation}\label{eq:P}
\min_{\beta\in\mathbb{R}^n,\,y\in\mathbb{R}^N} \left\{ \|y\| + \lambda p(\beta) \,|\, X\beta-Y=y \right\}.
\end{equation}
The dual problem of \eqref{eq:P} is given by
\begin{equation}
	- \min_{u\in\mathbb{R}^N} \left\{\inner{Y,u}+ (\lambda p)^*(X^Tu) + \delta_{\mathcal{B} }(u)\right\} \label{eq:D}
\end{equation}
where $\mathcal{B}$ is the unit Euclidean ball.

Compared with the Lasso-type problem \eqref{eq:lassotype}, the square-root regularized problem is more challenging to solve since both the loss function $\|\cdot\|$ and the regularizer $p$ are nonsmooth. We aim to use the framework of a proximal point algorithm (PPA) \citep{rockafellar1976monotone} for solving \eqref{eq:P}.
Given two sequences of positive parameters $\{\sigma_k\}$ and $\{\tau_k\} $ such that $\sigma_k\downarrow  \underline{\sigma}>0 $ and $\tau_k\downarrow \underline{\tau}>0 $, and an initial point $\left(\beta^0,y^0\right)\in\mathbb{R}^n \times \mathbb{R}^N$, the PPA for solving \eqref{eq:P} generates a sequence $\left\{\left(\beta^{k+1},y^{k+1}\right)\right\}$ via
\begin{equation}\label{eq:subPPA}
(\beta^{k+1},y^{k+1}) \approx \mathcal{P}_k(\beta^k,y^k):= \argmin{\beta\in\mathbb{R}^n,\,y\in\mathbb{R}^N}   \left\{
\begin{array}{c}
\|y\| + \lambda p({\beta})
\\[5pt]
+ \frac{\sigma_k}{2}\|\beta-\beta^k\|^2 +\frac{\tau_k}{2}\|y-y^k\|^2
\end{array}
\,\Bigg|\, X\beta -Y = y \right\}.
\end{equation}

\begin{algorithm}[H] \label{algo:PPA}
	\caption{Proximal point algorithm for solving \eqref{eq:P}}
	
	\textbf{Input} Data $X\in\mathbb{R}^{N\times n}$, $Y\in\mathbb{R}^{N}$, a penalty parameter $\lambda>0$, and a regularizer $p$.
	
	\textbf{Initialize} $\beta^0 = \bm{0}\in\mathbb{R}^n$, $y^0={-Y}$, $\sigma_0=\tau_0=1$.
	
	\While{the termination criterion is not met,}{ 		
		\textbf{Step 1.} Update $\beta^{k+1} $ and $y^{k+1}$ by solving  \eqref{eq:subPPA}.
		
		\textbf{Step 2.} Update $\sigma_{k+1}$ and $ \tau_{k+1} $. Set $k \leftarrow k+1 $.}
	
	\textbf{Output} $\beta^{k}$ and $y^k$.
\end{algorithm}
We use the standard criterion by \citet{rockafellar1976monotone} for controlling the inexactness
when solving \eqref{eq:subPPA} in Algorithm~\ref{algo:PPA}:
\begin{equation*}
\|(\beta^{k+1},y^{k+1}) - \mathcal{P}_k(\beta^k,y^k)\| \leq \varepsilon_k,\,\,\varepsilon_k > 0,\,\,\sum_{k\geq 0} \varepsilon_k < +\infty.
\end{equation*}
The global convergence of Algorihm~\ref{algo:PPA} follows from \citep{rockafellar1976monotone} directly.
A key difficulty in Algorithm~\ref{algo:PPA} is how to solve \eqref{eq:subPPA} efficiently.
Given $\sigma_k >0$, $\tau_k>0$, $\beta^k \in\mathbb{R}^n$, and $y^k\in\mathbb{R}^N$,
we recall the subproblem \eqref{eq:subPPA} given by
\begin{equation}\label{eq:subPPA2}
\begin{array}{cl}
\displaystyle \min_{\beta,y} & \displaystyle \|y\| + \lambda p({\beta}) + \frac{\sigma_k}{2}\|\beta-\beta^k\|^2 +\frac{\tau_k}{2}\|y-y^k\|^2  \\[6pt]
{\rm s.t.} & X\beta -Y = y.
\end{array}
\end{equation}
The Lagrangian function associated with  \eqref{eq:subPPA2} is given by
\begin{equation*}	
\begin{array}{cl}
\displaystyle \mathcal{L}(\beta,y;u) 		
 &=  \displaystyle \norm{y} + \frac{\tau_k}{2}\norm{y-y^k-\tau_k^{-1}u}^2 - \inner{y^k,u}-\frac{1}{2\tau_k}\norm{u}^2 - \inner{Y,u} + \lambda p({\beta})\\
\displaystyle &\quad + \frac{\sigma_k}{2}\norm{\beta-\beta^k+\sigma_k^{-1}X^Tu}^2 + \inner{\beta^k,X^Tu} - \frac{1}{2\sigma_k}\norm{X^Tu}^2,\,\,
\end{array}
\end{equation*}
for $(\beta,y,u)\in\mathbb{R}^n\times\mathbb{R}^N\times\mathbb{R}^N.$
By some simple manipulations, we can obtain that  the dual problem of \eqref{eq:subPPA2}, i.e., $\max_u \min_{\beta,y} \mathcal{L}(\beta,y;u)$, is given by
\begin{equation} \label{eq:subPPA-D}
\begin{array}{cl}
\displaystyle -\min_{u\in\mathbb{R}^N} &\displaystyle
\Psi(u):=\frac{1}{2\tau_k}\|u\|^2
+ \frac{1}{2\sigma_k}\|X^Tu\|^2 -\langle u,X\beta^k - y^k - Y\rangle\\[8pt]
&\displaystyle ~~~~~~~~~~
-\sigma_k \mathtt{M}_{\frac{\lambda}{\sigma_k} p}\left(\beta^k-\sigma_k^{-1}X^Tu \right) - \tau_k \mathtt{M}_{\frac{1}{\tau_k}\|\cdot\|} \left(y^k+\tau_k^{-1}u\right).
\end{array}
\end{equation}
It follows from the properties of Moreau envelope that the dual objective function $\Psi$ is continuously differentiable and convex.
Moreover, the gradient of $\Psi$ is given by
\begin{equation} \label{eq:gradient}
\nabla\Psi(u) = -X \mathtt{prox}_{\frac{\lambda}{\sigma_k} p}\left(\beta^k-\sigma_k^{-1}X^Tu \right)+ \mathtt{prox}_{\frac{1}{\tau_k}\|\cdot\|} \left(y^k+\tau_k^{-1}u\right)+ Y
,\,\,u\in\mathbb{R}^N.
\end{equation}
Thus the problem \eqref{eq:subPPA-D} is an unconstrained smooth convex minimization problem.
We let a dual optimal solution be $\bar{u}\in\arg\min \Psi(u)$.
Then the optimal solution $(\bar{\beta},\bar{y})$ to \eqref{eq:subPPA2} can be computed by
\begin{equation*}
\bar{\beta} = \mathtt{prox}_{\frac{\lambda}{\sigma_k} p}\left(\beta^k-\sigma_k^{-1}X^T\bar{u} \right),\quad
\bar{y}= \mathtt{prox}_{\frac{1}{\tau_k}\|\cdot\|} \left(y^k+\tau_k^{-1}\bar{u}\right).
\end{equation*}

Due to the favourable property of the dual problem \eqref{eq:subPPA-D}, we propose to solve \eqref{eq:subPPA2} via its dual.
In particular, the optimal solution of \eqref{eq:subPPA-D} is nothing but the solution of the nonlinear system $\nabla \Psi(u)=0,\,u\in\mathbb{R}^N$.
The latter can be solved by a semismooth Newton (SSN) method.
In order to apply the SSN method, we have to characterize
a certain generalized Jacobian of $\nabla\Psi(\cdot)$, which in turns depends on the generalized
Jacobian of $\mathtt{prox}_{\frac{\lambda}{\sigma_k}p}(\cdot)$.
Since the proximal mappings $\mathtt{prox}_{\frac{\lambda}{\sigma_k} p}\left(\cdot \right)$ and $\mathtt{prox}_{\frac{1}{\tau_k}\|\cdot\|} \left(\cdot\right)$ are Lipschitz continuous, the following multifunction, which is considered as a generalized Jacobian of
$\nabla\Psi(u)$, is well defined:
\begin{eqnarray*}
 \widehat{\partial} (\nabla \Psi)(u)
 = \left\{\sigma_k^{-1} X  U X^T + \tau_k^{-1} V \,\Bigg|
 \begin{array}{l}
U\in \partial\mathtt{prox}_{\frac{\lambda}{\sigma_k} p}\left(\beta^k-\sigma_k^{-1}X^Tu \right),
\\
V\in \partial\mathtt{prox}_{\frac{1}{\tau_k}\|\cdot\|} \left(y^k+\tau_k^{-1}u\right)
\end{array}
\right\}.
\end{eqnarray*}
Once an element $H \in \widehat{\partial} ( \nabla \Psi)(u) $ can be constructed explicitly for any given $u$, the SSN method can be implemented as follows.

\begin{algorithm}[H] \label{algo:SSN}
	\caption{Semismooth Newton method for solving \eqref{eq:subPPA-D}}
	
	\textbf{Input} Data $X\in\mathbb{R}^{N\times n}$, $Y\in\mathbb{R}^{N}$, a penalty parameter $\lambda>0$,  a regularizer $p$. $\beta^k\in\mathbb{R}^n$, $y^k\in\mathbb{R}^N$, $\sigma_k>0$, $\tau_k>0$. $\eta\in(0,1)$, $\varrho \in (0,1]$.  $\rho \in (0,1)$, $\mu\in(0,0.5)$.
 $\mathtt{tol} > 0$.
	
	\textbf{Initialize} $u^0=\bm{0}\in\mathbb{R}^{N}$, {$ f^0 = \nabla \Psi(u^{0})$ by \eqref{eq:gradient}}, $j=0$.
	
	\While{$\norm{{f^j}} > \mathtt{tol}$,}{ 		
		\textbf{Step 1.} Find an element ${H}^j \in \widehat{\partial} ( \nabla \Psi)(u^j) $, and then find an approximate solution $d^j$ to the linear system
		\begin{equation*}
		{H}^j d = -{f^j}
		\end{equation*}
such that $\|{H}^j d +{f^j}\| \leq \min (\eta,\| {f^j} \|^{1+\varrho})$. 		

		\textbf{Step 2.} Find a step size $\alpha_j = \rho^{m_j}$, where $m_j$ is the smallest nonnegative integer $m$ for which
		\begin{equation*}
		\Psi (u^j+\rho^{m_j} d^j) \leq \Psi(u^j) + \mu \rho^{m_j} \langle\nabla \Psi (u^j),d^j\rangle.
		\end{equation*}
		
		\textbf{Step 3.} Update
		$
        u^{j+1} = u^j +\alpha_j d^j
        $
        {and compute $f^{j+1} =  \nabla \Psi(u^j)$ by \eqref{eq:gradient}}.

        \textbf{Step 4.} Set $j \gets j+1$.
	}
	
	\textbf{Output} $u^j$.
\end{algorithm}

%
%
%
%
%
%

\bigskip

We can show that if
the optimal solution $(\bar{\beta},\bar{y})$ to \eqref{eq:subPPA2} does not overfit, i.e, $\bar{y} = X\bar{\beta} - Y \neq \bm{0}$,
then at the optimal solution $\bar{u}$ to \eqref{eq:subPPA-D},
the generalized Jacobian $\partial\mathtt{prox}_{\frac{1}{\tau_k}\|\cdot\|} \left(y^k+\tau_k^{-1}\bar{u}\right)$ is a singleton and
the element is positive definite.
But as the proof follows a similar argument to the one
in \citep[Proposition~12]{tang2020sparse}, we omit it here.
The above property is crucial to guarantee the fast convergence of the SSN method for solving \eqref{eq:subPPA-D}. We state the standard convergence result of the SSN method (Algorithm~\ref{algo:SSN}) without proof.

\begin{theorem}
Let $(\bar{\beta},\bar{y})$ be the optimal solution to the problem \eqref{eq:subPPA2}. Assume that the optimal solution does not overfit the data, i.e., $\bar{y} = X\bar{\beta} - Y \neq \bm{0}$. Then the sequence $\{u^j\}$ generated by Algorithm~\ref{algo:SSN} converges globally to the unique solution $\bar{u}$ of \eqref{eq:subPPA-D}. Furthermore, the local rate of convergence is of order $1+\varrho$, with $\varrho \in (0,1]$ given in Algorithm~\ref{algo:SSN}, i.e., for all $j$ sufficiently large,
$
\|u^{j+1} - \bar{u}\| = \mathcal{O}(\|u^j - \bar{u}\|^{1+\varrho}).
$
\end{theorem}

From now on, we restrict our discussions  to the case where the regularizer $p$ is either the sparse group Lasso regularizer \eqref{group-reg} or the fused Lasso regularizer \eqref{fused-reg}.
We will illustrate the explicit form of an element
$H \in \widehat{\partial} ( \nabla \Psi)(\cdot) $  based on \citep{li2018efficiently,zhang2020efficient}.
{It is worth noting that the matrix $H$ given below has (structured) sparsity inherited from that of matrices in the generalized Jacobian $\partial\mathtt{prox}_{p}\left( \cdot \right)$, and it is known as second order sparsity. The second order sparsity will reduce substantially the computational cost in the SSN method, as demonstrated in \citep{li2018highly,li2018efficiently,lin2019efficient,luo2019solving,zhang2020efficient}.}

{\bf Sparse group Lasso regularizer.}
We will first illustrate  the construction of an element in
$\widehat{\partial} ( \nabla \Psi)(\cdot)$ when $p$ is the sparse group Lasso regularizer given by \eqref{group-reg}.
Given positive constants $\kappa_1,\kappa_2$ and $\beta\in\mathbb{R}^n$, it is well known that for any $\beta$
\begin{equation*}
\mathtt{prox}_{\kappa_1 \|\cdot\|_1} (\beta) =  \left(|\beta| - \kappa_1 \bm{1} \right)^+ \odot \mathtt{sign}(\beta), \qquad
\mathtt{prox}_{\kappa_2 \|\cdot\|} (\beta) =
\begin{cases}
\left( 1 - \frac{\kappa_2}{\|\beta\|}\right)^+ \beta, & \mbox{if } \beta\neq 0, \\
\bm{0}, & \mbox{if } \beta = 0.
\end{cases}
\end{equation*}
We can construct a matrix $\bm{U}_{\kappa_1}(\beta) $ in $ \partial \mathtt{prox}_{\kappa_1 \|\cdot\|_1}(\beta)$ and a matrix $\bm{V}_{\kappa_2}(\beta)$ in $\partial \mathtt{prox}_{\kappa_2 \|\cdot\|}(\beta)$ for any $\beta$ respectively as follows
\begin{equation*}
\begin{array}{cl}
\partial \mathtt{prox}_{\kappa_1 \|\cdot\|_1} (\beta) &\ni \bm{U}_{\kappa_1}(\beta):=  \mathtt{Diag}(\bm{v}),
\\[6pt]
\partial \mathtt{prox}_{\kappa_2 \|\cdot\|} (\beta) &\ni \bm{V}_{\kappa_2}(\beta) :=
\begin{cases}
\left(1 - \frac{\kappa_2}{\|\beta\|} \right)I + \frac{\kappa_2}{\|\beta\|^3} \beta\beta^T, & \mbox{if }  \|\beta\| > \kappa_2,\\
\bm{0}, & \mbox{if }  \|\beta\| \leq \kappa_2,
\end{cases}
\end{array}
\end{equation*}
where $\bm{v}\in \mathbb{R}^n$ is defined by
$\bm{v}_i = 1\mbox{ if } |\beta_i| > \kappa_1$, and $\bm{v}_i=0$ otherwise.
Let $\tilde{\beta} :=\beta^k-\sigma_k^{-1}X^Tu  $ and $\tilde{y}:=y^k+\tau_k^{-1}u $. By \citep[(10) and Theorem~3.1]{zhang2020efficient},  we can construct a matrix $H \in \widehat{\partial} (\nabla\Psi)(\cdot)$ as follows
\begin{equation*}
	H := \sigma_k^{-1} \sum_{l=1}^g \left[ X_{G_l} \bm{V}_{\frac{\lambda w_2\omega_l}{\sigma_k}}\Big(  \mathtt{prox}_{\frac{\lambda w_1}{\sigma_k} \|\cdot\|_1} ({\tilde{\beta}_{G_l}}) \Big)  \bm{U}_{\frac{\lambda w_1}{\sigma_k} }({\tilde{\beta}_{G_l}} )  X^T_{G_l}\right]  +\tau_k^{-1} \bm{V}_{\frac{1}{\tau_k}}(\tilde{y}).
\end{equation*}

{\bf Fused Lasso regularizer.} Next, we  illustrate  the construction of an element in $\widehat{\partial} ( \nabla \Psi)(\cdot)$ when $p$ is the fused Lasso regularizer given by \eqref{fused-reg}.
From \citep{li2018efficiently}, we can construct a matrix $\bm{W}_{\kappa_2}(\beta)$ in
$\partial \mathtt{prox}_{\kappa_2 \norm{B\cdot}_1}(\beta) $ for any $\beta\in\mathbb{R}^n$ as follows
\begin{equation*}
	\partial \mathtt{prox}_{\kappa_2 \norm{B\cdot}_1}(\beta) \ni \bm{W}_{\kappa_2}(\beta):= I - B^T(\bm{\Sigma} B B^T \bm{\Sigma})^{\dagger}B,
\end{equation*}
where $\bm{\Sigma} = \mathtt{Diag}(\bm{\sigma}) $ and $\bm{\sigma}$ is given by
\begin{equation*}
	\bm{\sigma}_i  =
\begin{cases*}
		1, & if $ \left(B \mathtt{prox}_{\kappa_2 \norm{B\cdot}_1}(\beta) \right)_i =0 $,\\
		0, & otherwise,
\end{cases*}
\quad  {i =  1,2,\dots,n-1.}
\end{equation*}
Let $\tilde{\beta} :=\beta^k-\sigma_k^{-1}X^Tu  $ and $\tilde{y}:=y^k+\tau_k^{-1}u $. By \citep[(22) and Theorem~2]{li2018efficiently},  we can construct a matrix $H \in \widehat{\partial} (\nabla\Psi)(\cdot)$ as follows

\begin{equation*}
	H := \sigma_k^{-1}X\bm{U}_{\frac{\lambda w_1}{\sigma_k}}\Big(  \mathtt{prox}_{\frac{\lambda w_2}{\sigma_k} \|B\cdot\|_1} ({\tilde{\beta}}) \Big)  \bm{W}_{\frac{\lambda w_2 }{\sigma_k} }({\tilde{\beta}} )X^T +\tau_k^{-1} \bm{V}_{\frac{1}{\tau_k}}(\tilde{y}).
\end{equation*}

\section{Numerical Experiments} \label{sec:numerical-result}
In this section, we aim to evaluate the performance of our proximal point dual semismooth Newton algorithm (PPDNA) for solving the square-root regularized problem \eqref{eq:square-root}, when the regularizer $p$ is chosen to be the sparse group Lasso regularizer \eqref{group-reg} or the fused Lasso regularizer \eqref{fused-reg}.

\subsection{Setup of the Experiments}
Let $\mathtt{tol}$ be the tolerance, $\mathtt{maxiter}$ be the maximum {iteration number}, and $\mathtt{maxtime}$ be the maximum running time. We terminate a method at the $k{\rm th}$ iteration if the $k{\rm th}$ iterative point $\beta^k$ satisfies {one of the following conditions:}

\begin{itemize}
	\item $\norm{X\beta^k-Y}\ne0,\, \bar{\beta}^k:=X^T \frac{X\beta^k-Y}{\norm{X\beta^k-Y}}, $ and
	\begin{equation}\label{eq:error-kkt}
		\Delta_{\rm kkt}^k :=
		\frac{\norm{ \beta^k - \mathtt{prox}_{\lambda p} (\beta^k - \bar{\beta}^k)  }}{1+\norm{\beta^k}+\norm{\bar{\beta}^k} }<\mathtt{tol};
	\end{equation}
	\item  $\norm{X\beta^k-Y}=0$; {in this case, if the primal objective value in \eqref{eq:P} (denoted as ${\rm obj_\eqref{eq:P}}$) and the dual objective value in \eqref{eq:D} (denoted as ${\rm obj_\eqref{eq:D}}$) are available, we report the relative duality gap
		\begin{equation}\label{eq:error-gap}
			\Delta_{\rm pd.gap}^k := \dfrac{\rm obj_\eqref{eq:P} - obj_\eqref{eq:D} }{1 + \abs{\rm obj_\eqref{eq:P}} + \abs{\rm obj_\eqref{eq:D}} };
		\end{equation}
		otherwise, we report the relative successive change
		\begin{equation}\label{eq:error-change}
			\Delta_{\rm var.gap}^k := \dfrac{\norm{\beta^k-\beta^{k-1} }}{1+\norm{\beta^k}+\norm{{\beta}^{k-1}}};
		\end{equation}
	}
	
	\item $ k>\mathtt{maxiter}$, or the total running time exceeds $\mathtt{maxtime}$.
\end{itemize}
For all the algorithms, we set $\mathtt{tol}$ to be $10^{-7}$ and $\mathtt{maxtime}$ to be 30 minutes. In addition, we set  $\mathtt{maxiter}$ for our algorithm to be $10^2$, and for other algorithms to be $10^6$. All the experiments are performed in {\sc MATLAB} (version 9.7) on a {\sc Windows} workstation (6-core, Intel Core i7-8750H @ 2.20GHz, 8 Gigabytes of	RAM).

For the data matrix $X\in\mathbb{R}^{N\times n}$, we require  each column of $X$ to be nonzero, i.e., $\sum_{i=1}^{N} X_{ij}^2 >0,\,j=1,2,\dots,n$. For all the data matrices in our experiments, we normalize the columns such that the diagonal entries of the matrix $ \frac{1}{N}X^TX $ are equal to one. That is, we  let $\bm{d}\in\mathbb{R}^n $ be defined by $\bm{d}_j = \sqrt{   \frac{N}{\sum_{i=1}^{N} X_{ij}^2 } },\,j=1,2,\dots,n$ and then normalize $X$ by
$X \leftarrow X\mathtt{Diag}(\bm{d})$. Such a normalization of the data has been considered in  \citep{bunea2013group,stucky2017sharp,blanchet2019robust,MR4159660}. For the sparse group Lasso regularizer, we always choose the weights of groups as $\omega_l = {\sqrt{|{G_l}|}}$.

We also give the following explanations for the entries in the tables of numerical results.
We report an estimation of  the number of nonzero elements in a  computed vector $\beta\in\mathbb{R}^n$ as follows
\begin{equation}\label{eq:nnz-cal}
	{\rm nnz}(\beta) := \argmin{1\leq j\leq n} \Big\{ {j}\,\Big| \,\sum_{i=1}^j \abs{\beta_i^{\downarrow}} \geq 0.999\norm{\beta}_1 \Big\},
\end{equation}
where $\beta^{\downarrow}$ is a vector obtained from $\beta$  by sorting its entries in nonincreasing order of magnitude. When $\beta$ has a group structure $\{G_l \}_{l=1}^{g} $, which forms a partition of $\{1,2,\dots,n\} $, {we define $b\in\mathbb{R}^g $ by} $b_l := \norm{\beta_{G_l}}, \,l = 1,\dots,g$ and report an estimation of  the number of nonzero groups of $\beta$ by applying \eqref{eq:nnz-cal} to the vector $b$:
${\rm nnzgrp}(\beta) := {\rm nnz}(b)$.
Besides, we denote \({\rm nnzB(\beta) } := {\rm nnz} (B\beta)\) where \(B\beta = \left(\beta_1-\beta_2,\beta_2-\beta_3,\dots,\beta_{n-1} - \beta_{n} \right)^T  \).
We show the number of outer PPA iterations in Algorithm~\ref{algo:PPA} and the total number of inner SSN iterations in Algorithm~\ref{algo:SSN} of our PPDNA in the format of ``outer iteration $|$ inner iteration'' under the iteration column. The running time is in the format of ``minutes:seconds''.  An entry ``00'' under the column ``time'' means that the computational time is less than 0.5 second. Based on the value of $\|X\beta^k-Y\|$,  we report under the error column $\Delta_{\rm kkt}$, $*\Delta_{\rm pd.gap}$, or $\#\Delta_{\rm var.gap} $, given by \eqref{eq:error-kkt}, \eqref{eq:error-gap}, and \eqref{eq:error-change}, respectively.

\subsection{Alternating Direction Method of Multipliers for Solving the Square-root Regularized Problem}

To justify the necessity of our second order based method PPDNA in Section~\ref{sec:algorithm}, here we develop two first order methods based on the highly popular
alternating direction method of multipliers (ADMM) \citep{glowinski1975approximation,gabay1976dual} framework to compare our PPDNA against them.

Now we describe the implementation of ADMM  for solving {an equivalent form of the square-root regularized problem \eqref{eq:P} and its dual problem \eqref{eq:D}. By introducing slack variables $\alpha\in\mathbb{R}^n$, $\xi\in\mathbb{R}^n$, and $x\in\mathbb{R}^N$,  we obtain their equivalent forms respectively as follows}
\begin{eqnarray}
	&&	\min_{\beta\in\mathbb{R}^n,\,y\in\mathbb{R}^N,\,\alpha\in\mathbb{R}^n}  \left\{\norm{y} + \lambda p({\alpha})\mid X\beta-Y-y=0,\beta-\alpha=0 \right\} ,
	\label{eq:P2}
	\\[5pt]
	&-&\min_{u\in\mathbb{R}^N,\,\xi\in\mathbb{R}^n,\,x\in\mathbb{R}^N} \left\{{\inner{Y,u} + (\lambda p)^*(\xi) +\delta_{\mathcal{B}}(x)} \mid -X^Tu-\xi = 0,u-x=0 \right\}.
	\label{eq:D2}
\end{eqnarray}

Given a positive scalar $\mu$, the augmented Lagrangian functions associated with \eqref{eq:P2} and \eqref{eq:D2} are respectively given by
\begin{align*}
	\mathcal{L}^{(\rm P)}_{\mu}(\beta,y,\alpha;u,\xi) &= \norm{y} + \lambda p({\alpha}) + \frac{\mu}{2} \norm{ X\beta-Y-y+\mu^{-1}u }^2 +
	\frac{\mu}{2}\norm{\beta-\alpha+\mu^{-1}\xi}^2
	\\
	&\quad
	- \frac{1}{2\mu}\norm{u}^2 - \frac{1}{2\mu}\norm{\xi}^2,
	\quad \forall \ (\beta,y,\alpha,u,\xi) \in\mathbb{R}^n\times \mathbb{R}^N\times \mathbb{R}^n\times \mathbb{R}^N\times\mathbb{R}^n, \\
	\mathcal{L}^{(\rm D)}_{\mu}(u,\xi,x;\beta,y) &= \inner{Y,u}  +\delta_{\mathcal{B}}(x) + (\lambda p)^*(\xi)  +\frac{\mu}{2}\norm{-X^Tu-\xi+\mu^{-1}\beta}^2+\frac{\mu}{2}\norm{u-x+\mu^{-1}y}^2
	\\
	&\quad  - \frac{1}{2\mu}\norm{\beta}^2 - \frac{1}{2\mu}\norm{y}^2,
	\quad \forall (u,\xi,x,\beta,y) \in\mathbb{R}^N\times \mathbb{R}^n\times \mathbb{R}^N\times \mathbb{R}^n\times\mathbb{R}^N.
\end{align*}
The ADMM for solving the primal problem \eqref{eq:P2} and the dual problem \eqref{eq:D2} is given respectively in Algorithm~\ref{algo:pADMM} and Algorithm~\ref{algo:dADMM}. For the convergence results, we refer the reader to \citep{chen2017efficient,fazel2013hankel}.
We make the following remarks on techniques for solving linear systems \eqref{eq:pADMM-linsys} and \eqref{eq:dADMM-linsys}, which is the most expensive part in ADMM frameworks. When $n<N$, we solve the $n\times n$ linear system \eqref{eq:pADMM-linsys} either by the Cholesky factorization or the preconditioned conjugate gradients method, depending on $n$. Otherwise, we apply the Sherman-Morrison-Woodbury formula \citep{MR3024913}
$$
(I+X^TX)^{-1} = I -X(I+XX^T)^{-1}X^T,
$$
and only compute the Cholesky factorization of a smaller $N\times N$ matrix $I+XX^T$.  For solving the $N\times N$ linear system \eqref{eq:dADMM-linsys} when $N<n$, we use either the Cholesky factorization or the preconditioned conjugate gradients method, depending on $N$. Otherwise, we apply the Sherman-Morrison-Woodbury formula
$$ (I+XX^T)^{-1} = I - X^T(I+X^TX)^{-1}X,$$
and compute the Cholesky factorization of $I+X^TX$. In any case, the Cholesky factorization (if needed) is only performed once at the beginning of the algorithm.

	\begin{algorithm}[h] \label{algo:pADMM}
		\caption{pADMM for solving  \eqref{eq:P2}}
		
		\textbf{Input} Data $X\in\mathbb{R}^{N\times n}$, $Y\in\mathbb{R}^{N}$, a penalty parameter $\lambda>0$,  a regularizer $p$.
		
		\textbf{Initialize} $k=0,y^0=\bm{0}\in\mathbb{R}^N,\alpha^0=\bm{0}\in\mathbb{R}^n,u^0=\bm{0}\in\mathbb{R}^N,\xi^0=\bm{0}\in\mathbb{R}^n; \rho=1.618,\mu>0 $.
		
		\While{the termination criterion is not met,}{  		
			\textbf{Step 1.} Compute $\beta^{k+1}$ by
			\begin{align}
				\beta^{k+1} &= \argmin{\beta} 	\frac{\mu}{2} \norm{ X\beta-Y-y^k+\mu^{-1}u^k }^2 + \frac{\mu}{2}\norm{\beta-\alpha^k+\mu^{-1}\xi^k}^2 \notag \\
				&= (I+X^TX)^{-1} \left[ X^T(Y+y^k-\mu^{-1}u^k) +(\alpha^k-\mu^{-1}\xi^k) \right]. \label{eq:pADMM-linsys}
			\end{align}
			\textbf{Step 2.} Compute $y^{k+1},\alpha^{k+1}$ by
			\begin{align*}
				y^{k+1} &=\argmin{y} \norm{y} +  \frac{\mu}{2} \norm{ X\beta^{k+1}-Y-y+\mu^{-1}u^k }^2\\
				&= \mathtt{prox}_{\mu^{-1}\norm{\cdot}}\left(X\beta^{k+1}-Y+\mu^{-1}u^k \right),\\
				\alpha^{k+1} &= \argmin{\alpha} \lambda p(\alpha) +\frac{\mu}{2}\norm{\beta^{k+1}-\alpha+\mu^{-1}\xi^k}^2
				= \mathtt{prox}_{\mu^{-1}\lambda p} \left(\beta^{k+1}+\mu^{-1}\xi^k \right).
			\end{align*}
			\textbf{Step 3.} Update $u^{k+1},\xi^{k+1}$ by
			\begin{align*}
				u^{k+1 }= u^k + \rho \mu (X\beta^{k+1}-Y-y^{k+1}),\qquad
				\xi^{k+1} = \xi^k +\rho\mu(\beta^{k+1}-\alpha^{k+1} ).
			\end{align*}		
			\textbf{Step 4.} Set $k \leftarrow k+1 $.
		}
		
		\textbf{Output} $\beta^{k}$.
	\end{algorithm}

	\begin{algorithm}[h] \label{algo:dADMM}
		\caption{dADMM for solving \eqref{eq:D2}}
		
		\textbf{Input} Data $X\in\mathbb{R}^{N\times n}$, $Y\in\mathbb{R}^{N}$, a penalty parameter $\lambda>0$,  a regularizer $p$.
		
		\textbf{Initialize} $k=0,u^0=\bm{0}\in\mathbb{R}^N,\xi^0=\bm{0}\in\mathbb{R}^n,x^0=\bm{0}\in\mathbb{R}^N,\beta^0=\bm{0}\in\mathbb{R}^n,y^0=\bm{0}\in\mathbb{R}^N; \rho=1.618,\mu>0$.
		
		\While{the termination criterion is not met,}{ 		
			\textbf{Step 1.} Compute $u^{k+1}$ by
			\begin{align}
				u^{k+1} &= \argmin{u}u^TY +\dfrac{\mu}{2}\norm{-X^Tu-\xi^k+\mu^{-1}\beta^k}^2 +\dfrac{\mu}{2}\norm{u-x^k+\mu^{-1}y^k}^2 \notag \\
				&= (I+XX^T)^{-1}\left[-\mu^{-1}Y +X(\mu^{-1}\beta^k-\xi^k)-(\mu^{-1}y^k-x^k) \right].	\label{eq:dADMM-linsys}
			\end{align}
			\textbf{Step 2.} Compute $\xi^{k+1},x^{k+1}$ by
			\begin{align*}
				\xi^{k+1} &= \argmin{\xi}  (\lambda p)^*(\xi)  +\dfrac{\mu}{2}\norm{-X^Tu^{k+1}-\xi+\mu^{-1}\beta^{k}}^2	\\
				&= \mathtt{prox}_{\mu^{-1}(\lambda p)^* } \left(-X^Tu^{k+1}+\mu^{-1}\beta^{k} \right)\\
				&= \left(-X^Tu^{k+1}+\mu^{-1}\beta^{k} \right) - \mu^{-1} \mathtt{prox}_{\mu\lambda p } \left(-\mu X^Tu^{k+1}+\beta^{k} \right),\\			
				x^{k+1} &=\argmin{x}\delta_{\mathcal{B}}(x) +\dfrac{\mu}{2}\norm{u^{k+1}-x+\mu^{-1}y^k}^2	
				=  \Pi_{\mathcal{B}}\left(u^{k+1}+\mu^{-1}y^k\right).
			\end{align*}
			\textbf{Step 3.} Update $\beta^{k+1},y^{k+1}$ by
			\begin{align*}
				\beta^{k+1 }= \beta^k + \rho \mu (-X^Tu^{k+1}-\xi^{k+1}), \qquad y^{k+1} = y^k +\rho\mu(u^{k+1}-x^{k+1}).
			\end{align*}		
			\textbf{Step 4.} Set $k \leftarrow k+1 $.
			
		}
		
		\textbf{Output} $\beta^{k}$.
	\end{algorithm}

\subsection{Comparison of Efficiency for Solving the Square-root Sparse Group Lasso Problem}
In this section, we conduct extensive experiments to demonstrate the efficiency of our PPDNA for solving the square-root sparse group Lasso problem where the regularizer $p$ is given by \eqref{group-reg}.  In particular, we compare our  PPDNA (Algorithm~\ref{algo:PPA}+Algorithm~\ref{algo:SSN}) with pADMM (Algorithm~\ref{algo:pADMM}), dADMM (Algorithm~\ref{algo:dADMM}), and
the S-TISP solver \citep{bunea2013group}. Note that the S-TISP solver is limited to solving the square-root group Lasso problem where the coefficients in the expression \eqref{group-reg} of the regularier $p$ can only be taken as $(w_1,w_2) = (0,1)$ or  $(w_1,w_2) = (1,0)$ in that solver. For all tables in this section, we denote PPDNA, pADMM, dADMM, and S-TISP  by ``PP'', ``pA'', ``dA'', and ``ST'', respectively.

\subsubsection{Synthetic Data} \label{sec:synthetic-group}
We first show the results on synthetic data sets following the data generation mechanism in \citep{yuan2006model,bunea2013group,blanchet2017distributionally}. We first choose a correlation matrix $\Sigma$ to be a Toeplitz matrix, i.e., $\Sigma_{ij} = 0.5^{|i-j|}$, the dimensions of which will be clear from the context. We construct {four} examples in our experiments. As can be seen from the following examples, the true regression vector $\beta_0$ designed in Example~3 is not only groupwise sparse, but also sparse within a group. In contrast, the nonzero groups for $\beta_0$ in Example~1 and 2 are dense. {Besides, based on Example~3, we construct Example~4a/4b for which the true regression vectors are denser.}

\noindent{\bf Example~1} \citep{bunea2013group} The dimensions of this example are set as  $N=1000$, $g=200$ or $2000$, $n=3g$. We assign every three adjacent entries to be in one group, i.e., $G_l=\{3l-2,3l-1,3l\},\,l=1,\dots,g$. The true regression vector $\beta_0$ only contains three nonzero groups, i.e., $(\beta_0)_{G_l} = (2.5,2.5,2.5)^T$ for $l=1,3,4$ and $(\beta_0)_{G_l} = \bm{0}$ otherwise. The predictor vectors $X_i\in\mathbb{R}^n,\,i=1,\dots,N$ are generated from the multivariate normal distribution $\mathcal{N}(\bm{0},\Sigma)$. The response variables $Y_i\in\mathbb{R},\,i=1,\dots,N$ are generated from the linear regression model \eqref{linear-reg} with $\epsilon_i \overset{\rm i.i.d.}{\sim} \mathcal{N}(0,1)$ and $\sigma=1$.

\noindent{\bf Example~2} \citep{yuan2006model,blanchet2017distributionally} The dimensions of this example are set as $N=500$ or $10000$, $g=160$, $n=3g$. The group structure is given by $G_l=\{l,l+g,l+2g\},\,l=1,\dots,g$. The true regression vector $\beta_0$ only contains two nonzero groups, i.e., $(\beta_0)_{G_3} = (1,1,1)^T$, $(\beta_0)_{G_6} = (2/3,-1,1/2)^T$, and $(\beta_0)_{G_l} = \bm{0},\,l\notin \{3,6\}$. We generate random vectors $Z_i\in\mathbb{R}^{g},\,i=1,\dots,N$ from $\mathcal{N}(\bm{0},\Sigma)$ and a random scalar $\omega\in\mathbb{R}$ from $\mathcal{N}(0,1)$.  Let $  A_i:=\frac{Z_i+\omega\bm{1}}{\sqrt{2}}\in\mathbb{R}^{g}$. Then the predictor vectors $X_i\in\mathbb{R}^n,\,i=1,\dots,N$  are chosen to be the concatenation of three vectors constructed from $ A_i$:  $X_i := (A_i^T, A_i^T\odot A_i^T,A_i^T\odot A_i^T\odot A_i^T )^T$.  The response variables $Y_i\in\mathbb{R},\,i=1,\dots,N$ are generated from the linear regression model \eqref{linear-reg} with $\epsilon_i \overset{\rm i.i.d.}{\sim} \mathcal{N}(0,1)$ and $\sigma=2$.

\noindent{\bf Example~3} The dimensions of this example are set as $N=4000$, $g=2000$, $n=3g$ or $N=500$, $g=3000$, $n=3g$. The true regression vector $\beta_0$ only contains four nonzero groups, i.e., $(\beta_0)_{G_3} = (1,0,1)^T$, $(\beta_0)_{G_6} = (2/3,-1,0)^T$, $(\beta_0)_{G_9} = (-1,0,-1/2)^T$, $(\beta_0)_{G_{12}} = (0,-1,0)^T$, and $(\beta_0)_{G_l} = \bm{0},\,l\notin \{3,6,9,12\}$.
The group structure $\{G_l\}$ and the data $\{X,Y\} $ are generated in the same way as in \textbf{Example 2}.

\noindent{\bf Example~4a/4b} The dimensions of this example are set as \(N=4000,g=2000,n=3g\). Let \(\beta_{0}  \) be the  true regression vector in \textbf{Example 3}.  For a positive integer $T$, we define $\beta^{(T)}\in\mathbb{R}^n$ as follows:
\begin{equation*}
    \beta^{(T)}_{G_l} = \begin{cases}
        \left(\beta_0\right)_{G_i} & \text{ if } l = 12k+i,\,\,k \in\{0,1,\dots,T-1\},\,\, i \in\{1,2,\dots,12\}, \\
        \bm{0} & \text{ otherwise}.
    \end{cases}
\end{equation*}
By this construction, we have \(\mathrm{nnz} \left( \beta^{(T)}\right) = 7T \) and \(\mathrm{nnzgrp} \left( \beta^{(T)}\right) = 4T \). In particular, we set \(\beta^{(10)}\) and \(\beta^{(100)}\) as the true regression vectors of \textbf{Example 4a} and \textbf{Example 4b}, respectively. The  group structure $\{G_l\}$ and the data set $\{ X,Y\} $ are generated  in the same way as in \textbf{Example 3}.

We choose \(\lambda \in \left\{\lambda_{\text{Bun}},\lambda_{\text{StG}},\lambda_{\text{BlG}} \right\}  \), where the three tuning parameters $\lambda_{\text{Bun}}$, $\lambda_{\text{StG}}$, and $\lambda_{\text{BlG}}$ are theoretically optimal values given in \citep{bunea2013group,stucky2017sharp,blanchet2019robust}. See Appendix~\ref{sec:theo-lambda} for the details of the choice of the parameters.
First, we set $w_1=0$ and $w_2=1$ so that we can compare with the S-TISP solver, and the comparison of PPDNA, pADMM, dADMM, and S-TISP on {Examples~1-4}  are reported in Table~\ref{tab:synthetic-4-solvers}.
We can see from Table~\ref{tab:synthetic-4-solvers} that our PPDNA successfully solves all instances {except for Example~4b} within several seconds. For Example~3 with $(N,n)=(4000,6000)$,
pADMM fails to achieve the prespecified accuracy within 30 minutes, dADMM {takes} approximately 20 minutes,  and the running time of S-TISP  {varies} from about 2 minutes to 12 minutes; however, our PPDNA merely takes 5 seconds.
In addition, we set $w_1=w_2=0.5$ and compare the three methods PPDNA, pADMM, and dADMM in Table~\ref{tab:synthetic-3-solvers}. It can be observed that our PPDNA outperforms both pADMM and dADMM for all instances in Table~\ref{tab:synthetic-3-solvers}; our PPDNA returns accurate solutions for all instances within one minute. Even though dADMM is faster than pADMM,
its running time grows significantly with the increase of dimensions.
{In addition, comparing Examples 3, 4a, 4b ($N=4000,n=6000,g=2000$) where the numbers of nonzeros of the true regression vectors are 7, 70, 700, respectively, we found that the performances of the algorithms  do not vary too much.}
By comparing Table~\ref{tab:synthetic-4-solvers} and Table~\ref{tab:synthetic-3-solvers}, one may find that the efficiency of our PPDNA remains steady for the group Lasso regularizer  (\(w_1=0 \)) and the sparse group Lasso regularizer (\( w_1 = 0.5 \)).

\begin{table}[htbp] \centering
	\fontsize{6pt}{10pt}
	\setlength{\tabcolsep}{3pt} \selectfont
	\begin{tabular}{|c|r|r|rrrr|rrrr|rrrr|}
		\hline
		\multicolumn{1}{|c|}{problem} & \multicolumn{1}{c|}{\(\lambda\)} & \multicolumn{1}{c|}{nnz$|$} & \multicolumn{4}{c|}{iter} & \multicolumn{4}{c|}{time} & \multicolumn{4}{c|}{error} \\
		\cline{4-15}
		\multicolumn{1}{|c|}{\Longunderstack{\((N,n)\); \(g\) }} & & \multicolumn{1}{c|}{nnzgrp} &
		PP & pA & dA & ST &
		PP & pA & dA & ST &
		PP & pA & dA & ST \\
		\hline
		\multirow{3}{*}{\rotatebox{0}{\Longunderstack{Example 1 \\ $(1000,600)$ \\ $200$}}}
		& 3.485 & 9$|$3 & 17$|$70 & 1489 & 896 & 46 & 01 & 02 & 01 & 00 & 3.4e-08 & 4.6e-08 & 4.5e-08 & 8.8e-08 \\
		& 9.193 & 9$|$3 & 13$|$41 & 1071 & 1291 & 123 & 00 & 01 & 02 & 00 & 7.7e-09 & 5.3e-08 & 1.2e-08 & 9.4e-08 \\
		& 4.929 & 9$|$3 & 15$|$66 & 1489 & 1079 & 43 & 01 & 02 & 02 & 00 & 4.4e-09 & 3.3e-08 & 1.9e-08 & 7.8e-08 \\
		\hline
		\multirow{3}{*}{\rotatebox{0}{\Longunderstack{Example 1 \\ $(1000,6000)$ \\ $2000$}}}
		& 3.929 & 9$|$3 & 13$|$28 & 14374 & 6235 & 84 & 04 & 03:27 & 56 & 00 & 3.0e-11 & 5.3e-08 & 2.1e-08 & 8.8e-08 \\
		& 9.837 & 9$|$3 & 11$|$28 & 11054 & 3211 & 63 & 04 & 02:45 & 30 & 00 & 2.0e-08 & 2.0e-08 & 3.9e-08 & 8.5e-08 \\
		& 5.605 & 9$|$3 & 11$|$24 & 27501 & 6058 & 77 & 04 & 06:42 & 54 & 00 & 5.2e-08 & 3.7e-08 & 3.3e-08 & 8.7e-08 \\
		\hline
		\multirow{3}{*}{\rotatebox{0}{\Longunderstack{Example 2 \\ $(500,480)$ \\ $160$}}}
		& 3.462 & 6$|$2 & 14$|$48 & 12501 & 2993 & 2178 & 00 & 09 & 02 & 00 & 4.1e-09 & 1.6e-08 & 9.9e-08 & 1.0e-07 \\
		& 9.193 & 6$|$2 & 13$|$43 & 6389 & 3041 & 529 & 00 & 04 & 02 & 00 & 6.1e-08 & 2.4e-08 & 1.9e-08 & 9.9e-08 \\
		& 4.601 & 6$|$2 & 14$|$45 & 10001 & 3421 & 1860 & 00 & 07 & 02 & 00 & 1.1e-09 & 1.8e-08 & 2.5e-08 & 1.0e-07 \\
		\hline
		\multirow{3}{*}{\rotatebox{0}{\Longunderstack{Example 2 \\ $(10000,480)$ \\ $160$}}}
		& 3.353 & 6$|$2 & 14$|$53 & 100001 & 9442 & 3527 & 08 & 09:49 & 01:59 & 14 & 4.5e-09 & 7.0e-08 & 1.0e-07 & 1.0e-07 \\
		& 9.020 & 6$|$2 & 14$|$52 & 60001 & 16444 & 2574 & 08 & 05:54 & 03:28 & 10 & 1.6e-09 & 9.9e-08 & 9.9e-08 & 1.0e-07 \\
		& 4.546 & 6$|$2 & 13$|$55 & 72501 & 10014 & 3053 & 08 & 07:09 & 02:12 & 12 & 3.6e-08 & 8.6e-08 & 9.9e-08 & 1.0e-07 \\
		\hline
		\multirow{3}{*}{\rotatebox{0}{\Longunderstack{Example 3 \\ $(500,9000)$ \\ $3000$}}}
		& 4.026 & 23$|$8 & 13$|$37 & 27501 & 9581 & 27607 & 05 & 04:50 & 01:00 & 02:14 & 2.8e-09 & 7.0e-08 & 3.7e-08 & 1.0e-07 \\
		& 10.008 & 3$|$1 & 12$|$28 & 95001 & 15001 & 5870 & 05 & 17:25 & 01:37 & 28 & 1.2e-09 & 1.7e-08 & 3.0e-08 & 1.0e-07 \\
		& 17.169 & 0$|$0 & 1$|$1 & 1 & 1 & 1 & 00 & 00 & 00 & 00 & 0.0e-00 & 0.0e-00 & 0.0e-00 & 0.0e-00 \\
		\hline
		\multirow{3}{*}{\rotatebox{0}{\Longunderstack{Example 3 \\ $(4000,6000)$ \\ $2000$}}}
		& 3.800 & 12$|$4 & 14$|$51 & 24123 & 19459 & 34664 & 05 & 30:00 & 16:50 & 12:51 & 1.7e-09 & 7.1e-05 & 7.4e-08 & 1.0e-07 \\
		& 9.760 & 12$|$4 & 14$|$48 & 23646 & 22501 & 21300 & 05 & 30:00 & 19:52 & 08:07 & 3.6e-11 & 4.8e-06 & 1.8e-08 & 1.0e-07 \\
		& 17.163 & 12$|$4 & 13$|$48 & 23978 & 25001 & 7288 & 05 & 30:00 & 21:48 & 02:54 & 7.6e-08 & 7.7e-06 & 6.5e-08 & 1.0e-07 \\
		\hline
		\multirow{3}{*}{\rotatebox{0}{\Longunderstack{Example 4a \\ $(4000,6000)$ \\ $2000$}}}
		& 4.005 & 111$|$41 & 21$|$118 & 34205 & 7409 & 36603 & 06 & 30:00 & 05:21 & 06:45 & 4.1e-09 & 1.8e-07 & 1.0e-07 & 1.0e-07 \\
		& 9.791 & 132$|$48 & 21$|$111 & 10001 & 19713 & 19748 & 03 & 08:49 & 13:57 & 03:35 & 1.3e-08 & 2.1e-08 & 3.2e-08 & 1.0e-07 \\
		& 16.910 & 110$|$37 & 20$|$90 & 22501 & 25001 & 18372 & 02 & 19:51 & 17:40 & 03:20 & 6.3e-09 & 4.1e-08 & 4.4e-08 & 1.0e-07 \\
		\hline
		\multirow{3}{*}{\rotatebox{0}{\Longunderstack{Example 4b \\ $(4000,6000)$ \\ $2000$}}}
		& 3.994 & 1256$|$447 & 18$|$103 & 33392 & 4233 & 29186 & 01:33 & 30:00 & 03:11 & 05:22 & 3.8e-09 & 4.9e-07 & 9.8e-08 & 1.0e-07 \\
		& 9.791 & 1142$|$387 & 21$|$127 & 10001 & 7271 & 29179 & 26 & 09:04 & 05:22 & 05:21 & 1.4e-08 & 7.5e-08 & 1.0e-07 & 1.0e-07 \\
		& 17.466 & 852$|$290 & 22$|$143 & 12501 & 13635 & 22256 & 20 & 11:13 & 09:51 & 04:04 & 1.6e-08 & 2.7e-08 & 2.3e-08 & 1.0e-07 \\
		\hline		
	\end{tabular}
	\caption{Square-root sparse group Lasso model on synthetic data sets {Examples~1-4} with \((w_1,w_2) = (0,1) \). }
	\label{tab:synthetic-4-solvers}
\begin{tabular}{|c|r|r|rrr|rrr|rrr|}
		\hline
		\multicolumn{1}{|c|}{{problem}} & \multicolumn{1}{c|}{\(\lambda\)} & \multicolumn{1}{c|}{{nnz$|$}} & \multicolumn{3}{c|}{{iter}} & \multicolumn{3}{c|}{{time}} & \multicolumn{3}{c|}{{error}} \\
		\cline{4-12}
		\multicolumn{1}{|c|}{\Longunderstack{\((N,n)\); \(g\) }} &  & \multicolumn{1}{c|}{{nnzgrp}} &
		PP & pA & dA & PP & pA & dA & PP & pA & dA  \\
		\hline
		\multirow{3}{*}{\rotatebox{0}{\Longunderstack{Example 1 \\ $(1000,600)$ \\ $200$}}}
		& 3.589 & 9$|$3 & 14$|$58 & 1489 & 896 & 01 & 02 & 02 & 1.6e-08 & 4.5e-08 & 4.5e-08 \\
		& 9.262 & 9$|$3 & 13$|$44 & 1071 & 1291 & 00 & 01 & 02 & 1.5e-08 & 4.9e-08 & 1.1e-08 \\
		& 5.018 & 9$|$3 & 18$|$79 & 1489 & 1081 & 01 & 02 & 02 & 2.4e-09 & 3.2e-08 & 1.8e-08 \\
		\hline
		\multirow{3}{*}{\rotatebox{0}{\Longunderstack{Example 1 \\ $(1000,6000)$ \\ $2000$}}}
		& 4.064 & 9$|$3 & 12$|$26 & 22501 & 6243 & 04 & 05:46 & 01:00 & 2.7e-10 & 8.2e-08 & 1.5e-08 \\
		& 9.906 & 9$|$3 & 11$|$28 & 9581 & 3211 & 04 & 02:28 & 31 & 4.2e-08 & 3.7e-08 & 8.9e-08 \\
		& 5.615 & 9$|$3 & 11$|$24 & 15001 & 6210 & 04 & 03:51 & 59 & 7.5e-08 & 2.6e-08 & 4.9e-09 \\
		\hline
		\multirow{3}{*}{\rotatebox{0}{\Longunderstack{Example 2 \\ $(500,480)$ \\ $160$}}}
		& 3.724 & 6$|$3 & 14$|$56 & 10001 & 2892 & 01 & 08 & 02 & 5.6e-09 & 1.2e-08 & 9.9e-08 \\
		& 9.298 & 6$|$2 & 14$|$51 & 5325 & 3062 & 01 & 04 & 02 & 3.2e-08 & 9.3e-08 & 1.2e-08 \\
		& 4.593 & 6$|$2 & 14$|$53 & 10001 & 3286 & 01 & 08 & 02 & 1.7e-09 & 1.4e-08 & 9.0e-08 \\
		\hline
		\multirow{3}{*}{\rotatebox{0}{\Longunderstack{Example 2 \\ $(10000,480)$ \\ $160$}}}
		& 3.454 & 6$|$3 & 18$|$81 & 52501 & 8454 & 10 & 05:29 & 01:53 & 1.9e-09 & 6.9e-08 & 1.0e-07 \\
		& 9.039 & 6$|$2 & 15$|$60 & 60001 & 13449 & 08 & 06:15 & 02:57 & 5.0e-10 & 9.4e-08 & 4.2e-08 \\
		& 4.564 & 6$|$2 & 16$|$79 & 72501 & 14916 & 11 & 07:33 & 03:22 & 6.1e-09 & 8.6e-08 & 9.9e-08 \\
		\hline
		\multirow{3}{*}{\rotatebox{0}{\Longunderstack{Example 3 \\ $(500,9000)$ \\ $3000$}}}
		& 4.422 & 22$|$9 & 13$|$34 & 32501 & 7986 & 05 & 05:49 & 54 & 2.1e-11 & 6.1e-08 & 6.3e-08 \\
		& 10.114 & 3$|$1 & 12$|$29 & 35001 & 11496 & 05 & 06:21 & 01:15 & 7.9e-10 & 6.2e-08 & 8.3e-08 \\
		& 17.216 & 0$|$0 & 1$|$1 & 1 & 1 & 00 & 00 & 00 & 0.0e-00 & 0.0e-00 & 0.0e-00 \\
		\hline
		\multirow{3}{*}{\rotatebox{0}{\Longunderstack{Example 3 \\ $(4000,6000)$ \\ $2000$}}}
		& 3.955 & 10$|$5 & 14$|$50 & 23062 & 22896 & 05 & 30:00 & 20:45 & 4.2e-10 & 8.0e-05 & 1.0e-07 \\
		& 9.791 & 10$|$4 & 14$|$54 & 22845 & 22046 & 05 & 30:00 & 20:05 & 1.3e-09 & 6.0e-06 & 3.6e-08 \\
		& 16.958 & 8$|$4 & 12$|$39 & 22889 & 22501 & 05 & 30:00 & 20:26 & 7.8e-08 & 2.7e-06 & 8.2e-08 \\
		\hline
		\multirow{3}{*}{\rotatebox{0}{\Longunderstack{Example 4a \\ $(4000,6000)$ \\ $2000$}}}
		& 4.005 & 115$|$72 & 20$|$106 & 34158 & 11104 & 03 & 30:00 & 08:09 & 3.5e-08 & 1.8e-06 & 1.0e-07 \\
		& 9.791 & 101$|$78 & 21$|$104 & 10001 & 19907 & 03 & 08:49 & 14:08 & 4.3e-08 & 5.6e-08 & 6.1e-08 \\
		& 16.910 & 133$|$61 & 21$|$98 & 15001 & 17501 & 03 & 13:14 & 12:23 & 1.5e-08 & 3.6e-08 & 4.2e-08 \\
		\hline
		\multirow{3}{*}{\rotatebox{0}{\Longunderstack{Example 4b \\ $(4000,6000)$ \\ $2000$}}}
		& 3.994 & 889$|$483 & 22$|$111 & 27501 & 4324 & 41 & 24:51 & 03:18 & 1.4e-08 & 8.2e-08 & 1.0e-07 \\
		& 9.791 & 1114$|$403 & 21$|$100 & 25001 & 5183 & 21 & 22:35 & 03:51 & 5.6e-09 & 9.2e-08 & 1.0e-07 \\
		& 17.466 & 884$|$298 & 18$|$130 & 7501 & 8190 & 19 & 06:45 & 05:56 & 8.6e-09 & 3.3e-08 & 6.9e-08 \\
		\hline
	\end{tabular}
	\caption{Square-root sparse group Lasso model on synthetic data sets {Examples~1-4} with \(w_1=w_2=0.5 \).}
	\label{tab:synthetic-3-solvers}
\end{table}

\subsubsection{UCI Data with Synthetic Group Structure}\label{sec:group:uci}
We use the UCI repository \citep{asuncion2007uci,chang2011libsvm} in this section.  Following \citep{huang2010predicting,li2018highly,zhang2020efficient}, we expand the original features of the data sets \textit{housing}, \textit{bodyfat}, \textit{pyrim}, \textit{triazines}. The UCI data sets are not equipped with group structures, and we design the group structure as follows. We set the total number of groups as $g=300$, and a feature is assigned randomly to a group with a uniform probability.

\begin{table}[!h] \centering
	\fontsize{6pt}{10pt}\setlength{\tabcolsep}{3pt} \selectfont
	\begin{tabular}{|r|r|r|rrrr|rrrr|rrrr|}
		\hline
		\multicolumn{1}{|c|}{\Longunderstack{problem}} & \multicolumn{1}{c|}{\(w_1|\lambda\)} & \multicolumn{1}{c|}{\Longunderstack{nnz$|$}} & \multicolumn{4}{c|}{\Longunderstack{iter}} & \multicolumn{4}{c|}{\Longunderstack{time}} & \multicolumn{4}{c|}{\Longunderstack{error}} \\
		\cline{4-15}
		\multicolumn{1}{|c|}{\Longunderstack{\((N,n);g\) }} & & \multicolumn{1}{c|}{\Longunderstack{nnzgrp}}
		& PP & pA & dA & ST & PP & pA & dA & ST & PP & pA & dA & ST  \\
		\hline
		
		\multirow{6}{*}{\rotatebox{0}{\Longunderstack{\textit{housing} \\ $(253,77520)$ \\ $300$}}} & 0.0$|$0.5 & 1846$|$8 & 12$|$32 & 3851 & 2470 & 86386 & 08 & 02:36 & 54 & 30:00 & 1.2e-08 & 3.8e-08 & 1.0e-07 & 1.7e-04 \\
		& 0.0$|$1.0 & 1640$|$7 & 17$|$66 & 7986 & 4341 & 88088 & 02 & 05:15 & 01:32 & 30:00 & 1.5e-08 & 4.1e-08 & 2.0e-08 & 8.8e-05 \\
		& 0.0$|$2.0 & 724$|$3 & 16$|$66 & 17501 & 7986 & 89932 & 03 & 11:26 & 02:45 & 30:00 & 1.9e-08 & 4.0e-08 & 4.2e-08 & 3.3e-05 \\
		\cline{2-15}
		& 0.5$|$0.5 & 1145$|$11 & 15$|$39 & 5546 & 2909 & $-$ & 09 & 03:43 & 01:06 & $-$ & 1.5e-10 & 2.6e-08 & 1.0e-07 & $-$ \\
		& 0.5$|$1.0 & 709$|$7 & 15$|$61 & 11496 & 4495 & $-$ & 02 & 07:39 & 01:37 & $-$ & 3.5e-09 & 5.5e-08 & 3.2e-08 & $-$ \\
		& 0.5$|$2.0 & 460$|$4 & 13$|$48 & 20001 & 9020 & $-$ & 01 & 12:57 & 03:07 & $-$ & 7.6e-09 & 7.5e-08 & 2.0e-08 & $-$ \\
		\hline
		\multirow{6}{*}{\rotatebox{0}{\Longunderstack{\textit{bodyfat} \\ $(126,116280)$ \\ $300$}}} & 0.0$|$0.5 & 703$|$2 & 12$|$25 & 11496 & 3081 & 104870 & 08 & 05:52 & 50 & 30:00 & 6.7e-11 & 9.9e-08 & 2.6e-08 & 6.2e-03 \\
		& 0.0$|$1.0 & 349$|$1 & 12$|$53 & 45001 & 17501 & 106770 & 01 & 22:54 & 04:42 & 30:00 & 1.4e-08 & 8.8e-08 & 5.4e-08 & 4.3e-03 \\
		& 0.0$|$2.0 & 348$|$1 & 12$|$39 & 60191 & 32501 & 107478 & 01 & 30:00 & 08:35 & 30:00 & 1.8e-11 & 2.2e-06 & 8.9e-08 & 3.3e-03 \\
		\cline{2-15}
		& 0.5$|$0.5 & 84$|$3 & 14$|$38 & 47501 & 10772 & $-$ & 08 & 25:02 & 03:21 & $-$ & 9.0e-11 & 9.7e-08 & 1.0e-07 & $-$ \\
		& 0.5$|$1.0 & 85$|$2 & 13$|$60 & 56682 & 8317 & $-$ & 02 & 30:00 & 02:24 & $-$ & 5.9e-08 & 8.3e-07 & 8.7e-08 & $-$ \\
		& 0.5$|$2.0 & 227$|$2 & 13$|$50 & 59609 & 27501 & $-$ & 01 & 30:00 & 07:20 & $-$ & 7.3e-09 & 1.0e-05 & 6.4e-08 & $-$ \\
		\hline
		\multirow{6}{*}{\rotatebox{0}{\Longunderstack{\textit{pyrim} \\ $(37,169911)$ \\ $300$}}} & 0.0$|$0.5 & 1072$|$2 & 13$|$27 & 45001 & 3211 & 129108 & 07 & 18:12 & 44 & 30:00 & 1.6e-09 & 6.5e-08 & 3.9e-08 & 1.3e-02 \\
		& 0.0$|$1.0 & 577$|$1 & 11$|$41 & 52501 & 7986 & 62705 & 01 & 21:16 & 01:46 & 14:05 & 3.1e-08 & 7.9e-08 & 1.1e-08 & 1.0e-07 \\
		& 0.0$|$2.0 & 1099$|$2 & 12$|$39 & 75001 & 13795 & 135011 & 01 & 29:58 & 03:03 & 30:00 & 1.9e-09 & 9.8e-08 & 1.6e-08 & 3.3e-03 \\
		\cline{2-15}
		& 0.5$|$0.5 & 653$|$3 & 11$|$21 & 37501 & 2300 & $-$ & 08 & 15:24 & 32 & $-$ & 3.1e-10 & 3.7e-08 & 5.6e-08 & $-$ \\
		& 0.5$|$1.0 & 215$|$1 & 12$|$48 & 32501 & 12501 & $-$ & 01 & 13:17 & 02:51 & $-$ & 9.0e-11 & 6.3e-08 & 8.6e-08 & $-$ \\
		& 0.5$|$2.0 & 350$|$1 & 11$|$38 & 67501 & 9581 & $-$ & 01 & 26:40 & 02:07 & $-$ & 8.1e-09 & 9.4e-08 & 8.2e-08 & $-$ \\
		\hline
		\multirow{6}{*}{\rotatebox{0}{\Longunderstack{\textit{triazines}\\ $(93,557845)$ \\ $300$}}}  & 0.0$|$0.5 & 7040$|$4 & 18$|$75 & 7961 & 6656 & 11474 & 20 & 30:00 & 12:39 & 30:00 & 9.8e-09 & 9.9e-06 & 5.2e-08 & 5.4e-02 \\
		& 0.0$|$1.0 & 3556$|$2 & 18$|$80 & 8000 & 14856 & 12325 & 15 & 30:00 & 30:00 & 30:00 & 7.2e-09 & 9.9e-06 & 1.1e-07 & 6.3e-02 \\
		& 0.0$|$2.0 & 1806$|$1 & 13$|$57 & 8981 & 17492 & 13194 & 10 & 30:00 & 30:00 & 30:00 & 1.1e-08 & 8.9e-06 & 1.7e-07 & 1.9e-02 \\
		\cline{2-15}
		& 0.5$|$0.5 & 2809$|$3 & 21$|$99 & 8630 & 5546 & $-$ & 20 & 30:00 & 09:56 & $-$ & 8.9e-08 & 3.7e-05 & 9.2e-08 & $-$ \\
		& 0.5$|$1.0 & 1987$|$2 & 21$|$94 & 8665 & 7986 & $-$ & 16 & 30:00 & 14:25 & $-$ & 2.9e-08 & 9.2e-06 & 3.0e-08 & $-$ \\
		& 0.5$|$2.0 & 2959$|$3 & 14$|$58 & 9016 & 17286 & $-$ & 10 & 30:00 & 30:00 & $-$ & 2.4e-09 & 8.8e-06 & 9.1e-08 & $-$ \\
		\hline
		\multirow{6}{*}{\rotatebox{0}{\Longunderstack{\textit{E2006.test} \\ $(1654,72812)$ \\ $300$}}} & 0.0$|$2.0 & 5972$|$37 & 15$|$69 & 37501 & 12501 & 17749 & 03 & 13:23 & 03:18 & 02:06 & 3.6e-08 & 6.0e-09 & 5.4e-08 & 9.7e-08 \\
		& 0.0$|$5.0 & 1395$|$8 & 14$|$46 & 11497 & 15001 & 4100 & 00 & 03:58 & 03:42 & 27 & 8.8e-11 & 7.6e-08 & 7.3e-08 & 9.9e-08 \\
		& 0.0$|$7.5 & 0$|$0 & 1$|$1 & 1 & 1 & 1 & 00 & 00 & 00 & 00 & 0.0e-00 & 0.0e-00 & 0.0e-00 & 0.0e-00 \\
		\cline{2-15}
		& 0.5$|$2.0 & 4641$|$65 & 19$|$82 & 32501 & 18610 & $-$ & 05 & 12:01 & 05:01 & $-$ & 5.6e-09 & 6.6e-08 & 4.6e-08 & $-$ \\
		& 0.5$|$5.0 & 656$|$22 & 13$|$49 & 35001 & 20692 & $-$ & 03 & 12:17 & 05:25 & $-$ & 1.3e-10 & 8.9e-08 & 1.1e-08 & $-$ \\
		& 0.5$|$7.5 & 149$|$6 & 12$|$42 & 15001 & 15001 & $-$ & 00 & 05:13 & 03:45 & $-$ & 6.5e-08 & 7.4e-08 & 8.1e-08 & $-$ \\
		\hline
		\multirow{6}{*}{\rotatebox{0}{\Longunderstack{\textit{E2006.train} \\ $(8044,150348)$ \\ $300$}}} & 0.0$|$2.0 & 16098$|$46 & 19$|$94 & $-$ & 3031 & 16517 & 05:31 & $-$ & 30:00 & 07:44 & 9.3e-10 & $-$ & 2.7e-01 & 9.6e-08 \\
		& 0.0$|$5.0 & 3088$|$9 & 17$|$71 & $-$ & 3010 & 11208 & 05 & $-$ & 30:00 & 05:03 & 2.2e-08 & $-$ & 4.6e-04 & 9.8e-08 \\
		& 0.0$|$7.5 & 1092$|$3 & 14$|$57 & $-$ & 3012 & 3967 & 02 & $-$ & 30:00 & 01:47 & 1.8e-09 & $-$ & 3.8e-04 & 1.0e-07 \\
		\cline{2-15}
		& 0.5$|$2.0 & 10427$|$73 & 20$|$100 & $-$ & 3015 & $-$ & 33 & $-$ & 30:00 & $-$ & 8.3e-08 & $-$ & 3.0e-01 & $-$ \\
		& 0.5$|$5.0 & 1070$|$23 & 17$|$71 & $-$ & 3006 & $-$ & 08 & $-$ & 30:00 & $-$ & 1.3e-08 & $-$ & 6.2e-04 & $-$ \\
		& 0.5$|$7.5 & 296$|$9 & 15$|$60 & $-$ & 2865 & $-$ & 04 & $-$ & 30:00 & $-$ & 9.0e-09 & $-$ & 2.1e-02 & $-$ \\
		\hline
	\end{tabular}
	
	\caption{Square-root sparse group Lasso model on UCI data sets with \((w_1,w_2) = (0,1) \) or \((0.5,0.5) \). ``$-$'' denotes that the method is not applicable for the instance. }
	\label{tab:UCI-4-solvers}
\end{table}

On each UCI data set, we consider two cases:  $(w_1,w_2)=(0,1)$ and $(w_1,w_2)=(0.5,0.5)$. For each case, we select three values of \(\lambda\) from the set \(\{0.5,1,2,5,7.5\} \)  so that the resulting solutions do not overfit. When \(w_1=0\), we compare the performances of PPDNA with pADMM, dADMM, and S-TISP. Otherwise, we compare the performances of PPDNA with pADMM and dADMM. The results are presented in Table~\ref{tab:UCI-4-solvers}. It can be seen that PPDNA is always the fastest among the compared methods for all instances in Table~\ref{tab:UCI-4-solvers}. We also find that the iteration numbers and computational time of pADMM, dADMM, and S-TISP  fluctuates wildly with the changes in the dimensions or parameters. In contrast, the performances of our PPDNA are generally robust; it solves all the problems in Table~\ref{tab:UCI-4-solvers} within 30 seconds except for \textit{E2006.train}. In addition,  we can observe that pADMM is generally slow for problems in Table~\ref{tab:UCI-4-solvers}, and it  can not solve the last problem due to  insufficient of memory for computing the Cholesky factorization. Compared to S-TISP, dADMM is faster for the first four problems. For the last two problems, dADMM fails to solve them within 30 minutes, which is due to the computational cost for solving the linear system \eqref{eq:dADMM-linsys}. But note that S-TISP does not solve some of instances to the required level of accuracy. Since the dimensions of problems in UCI data sets are much larger than those of the synthetic data sets in Section~\ref{sec:synthetic-group},  pADMM, dADMM, and S-TISP reach the maximum running time 30 minutes in many instances. In contrast, our PPDNA succeeds in solving all instances. We can safely conclude that our PPDNA can be more efficient than pADMM, dADMM, and S-TISP for solving large-scale square-root sparse group Lasso problems.

\subsubsection{Real Data} \label{sec:real-data-group}
In this section, we present the numerical results of the square-root sparse group Lasso model on two real data sets which are equipped with natural group structures.
For a given data set $(X,Y)$ in this section, we randomly split it into the training set $(X_{\rm train}, Y_{\rm train})$ and the test set $(X_{\rm test}, Y_{\rm test})$ so that the number of observations \(N_{\rm train}\)  of the training data set is roughly twice larger than the number of observations \(N_{\rm test}\) of the test set.
Based on the training set, we first set $(w_1,w_2) = (0,1)$ and conduct 8-fold cross validation (CV) for selecting $\lambda$ over the set
\begin{equation}\label{CV2}
	\lambda \in \{10^{-1},10^{-0.95},10^{-0.9},\dots, 10^{0.95},10^1  \}.
\end{equation}
We then conduct 8-fold CV for selecting parameters $w_1$, $w_2$, and $\lambda$ over the sets
\begin{equation}\label{CV1}
	w_1 \in \{0,0.1,0.2,\dots,0.9,1\},\, w_2 = 1-w_1, \text{ and }\lambda \in \{10^{-1},10^{-0.95},10^{-0.9},\dots, 10^{0.95},10^1  \}.
\end{equation}
Based on the test set, we report the mean squared error (MSE) for an approximate solution $\beta$ defined by
\({\norm{X_{\rm test}\beta - Y_{\rm test} }^2 }/{N_{\rm test }} \).

\bigskip

\noindent{\bf Climate data} \citep{kalnay1996ncep}
The data records climate information of 10512 locations across the globe ($73\times144 $ grid of latitude and longitude, resolution $2.5^0\times 2.5^0$). For each location, it records the monthly means of 7 predictor variables {\it Air Temperature}, {\it Precipitable water}, {\it Relative humidity}, {\it Pressure}, {\it Sea Level Pressure}, {\it Horizontal	Wind Speed}, and {\it Vertical Wind Speed}. The predictor vector $X_i\in\mathbb{R}^{73584}$ is the concatenation of the 7 predictor variables at 10512 locations in the $i$th month, and we have data from 1948/1/1 to 2018/5/31 containing $N=814$ months. We {regard} the 7 predictor variables at a location as a group. The response variable is designed as follows. We first select a target location and then set $Y_i$ to be the {\it Air Temperature} at the target location in the $i$th month. Moreover, we remove the 7 predictor variables corresponding to the selected target location in $X$ and eventually we have the data $X\in\mathbb{R}^{814 \times 73577}$ and $Y\in\mathbb{R}^{814}$.

\begin{figure}[htbp]
	\centering
	\includegraphics[width=\textwidth]{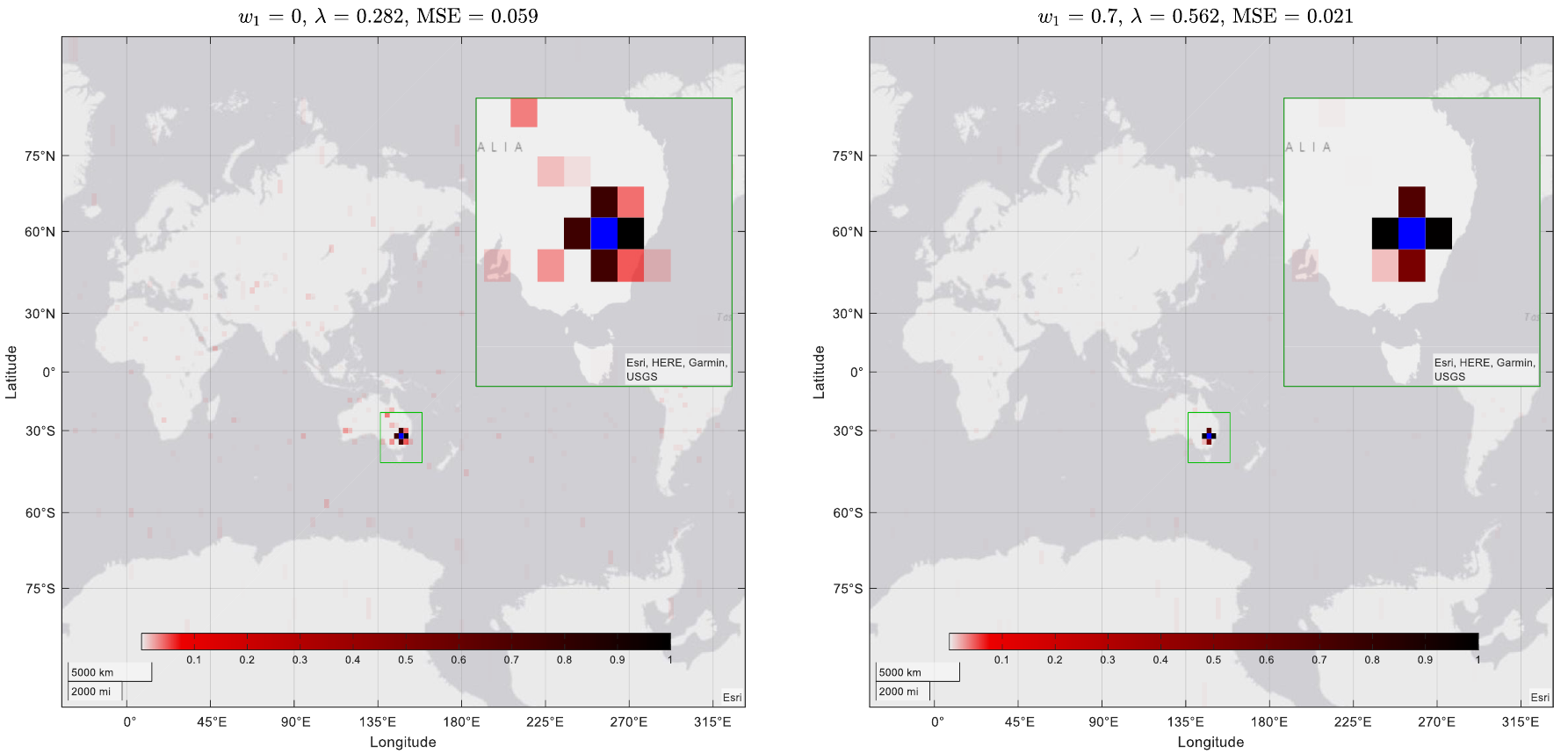}
	\caption{Active groups to predict \textit{Air Temperature} in a neighborhood of Dubbo (in blue).
		Left: \((w_1,w_2,\lambda)=(0,1,0.282) \).
		Right: \((w_1,w_2,\lambda)=(0.7,0.3,0.562) \). At the $l$th location, the value of ${\norm{\beta_{G_l}}}/{\max_{1\leq l \leq g}\norm{\beta_{G_l}}} $ is displayed.}
	\label{fig:climate}
\end{figure}

This data has also been used in \citep{ndiaye2016gap,zhang2020efficient}, showing that the sparse group Lasso regularizer is suitable for prediction in climate data. We consider five different target locations for the diversity of experiments. We first set $(w_1,w_2)=(0,1)$, and select $\lambda$ by CV over \eqref{CV2}. In addition, we select $w_1$, $w_2$, and $\lambda$ by CV over \eqref{CV1}. Table~\ref{tab:climate-4-solvers} shows the comparisons of PPDNA, pADMM, dADMM, and S-TISP on the climate data sets with five different target locations. As shown in Table~\ref{tab:climate-4-solvers}, our PPDNA outperforms the other three methods pADMM, dADMM, and S-TISP by a wide margin in term of computational time. In particular, both pADMM and S-TISP {fail} to solve any instance within 30 minutes, and the accuracy of the solutions returned by S-TISP, approximately $10^{-3}$,  are still far from being {satisfactory}. In addition, we plot the active groups predicting \textit{Air Temperature} in a neighborhood of Dubbo, New South Wales, Australia in Figure~\ref{fig:climate}.   As one can expect, the active groups shown in the Figure~\ref{fig:climate} are close or contiguous to the target location Dubbo.

\begin{table}[htbp] \centering
	\fontsize{6pt}{10pt}\setlength{\tabcolsep}{3pt} \selectfont
	\renewcommand{\arraystretch}{2}
	\begin{tabular}{|r|r|r|rrrr|rrrr|rrrr|}
		\hline
		\multicolumn{1}{|c|}{{target location}} & \multicolumn{1}{c|}{\(w_1|\lambda\)} &  \multicolumn{1}{c|}{MSE}  & \multicolumn{4}{c|}{{iter}} & \multicolumn{4}{c|}{{time}} & \multicolumn{4}{c|}{{error}} \\
		\cline{4-15}
		\multicolumn{1}{|c|}{{\((N_{\rm train},n), g\) }} & &  &    PP & pA & dA & ST & PP & pA & dA & ST &PP & pA & dA & ST \\
		\hline
		\multirow{2}{*}{\rotatebox{0}{\Longunderstack{Dakar \\\( (15^{\circ}N, 17^{\circ}30'W) \) \\ $(563,73577),\  10511$}}}  & 0.0$|$0.501 & 0.294 & 16$|$80 & 19469 & 6011 & 35469 & 17 & 30:00 & 05:56 & 30:00 & 1.3e-08 & 1.2e-05 & 9.8e-08 & 5.4e-03 \\
		\cline{2-15}
		& 1.0$|$1.259 & 0.028 & 19$|$109 & 20467 & 17760 & 32782 & 09 & 30:00 & 16:21 & 30:00 & 6.3e-08 & 8.5e-06 & 1.0e-07 & 1.8e-02 \\
		\hline
		\multirow{2}{*}{\rotatebox{0}{\Longunderstack{Dubbo \\\( (32^{\circ}30'S, 147^{\circ}30'E) \) \\ $(563,73577),\  10511$}}}& 0.0$|$0.282 & 0.059 & 16$|$45 & 19898 & 7245 & 35554 & 47 & 30:00 & 07:14 & 30:00 & 2.3e-08 & 1.5e-05 & 9.9e-08 & 3.8e-03 \\
		\cline{2-15}
		& 0.7$|$0.562 & 0.021 & 20$|$140 & 20235 & 11318 & $-$ & 13 & 30:00 & 11:09 & $-$ & 4.5e-08 & 5.4e-06 & 1.0e-07 & $-$ \\
		\hline
		\multirow{2}{*}{\rotatebox{0}{\Longunderstack{Enshi \\\( (30^{\circ}N, 110^{\circ}E) \) \\ $(563,73577),\  10511$}}} & 0.0$|$0.282 & 0.031 & 1$|$186 & 20230 & 7397 & 35177 & 58 & 30:00 & 07:34 & 30:00 & 6.2e-07 & 8.2e-06 & 1.0e-07 & 3.9e-03 \\
		\cline{2-15}
		& 0.8$|$0.794 & 0.023 & 21$|$123 & 20436 & 19163 & $-$ & 08 & 30:00 & 18:39 & $-$ & 5.3e-08 & 9.3e-06 & 1.0e-07 & $-$ \\
		\hline
		\multirow{2}{*}{\rotatebox{0}{\Longunderstack{Weihai \\\( (37^{\circ}30'N, 122^{\circ}30'E) \) \\ $(563,73577),\  10511$}}} & 0.0$|$0.282 & 0.045 & 19$|$74 & 20220 & 9667 & 31370 & 51 & 30:00 & 09:52 & 30:00 & 1.1e-10 & 3.6e-06 & 1.0e-07 & 4.1e-03 \\
		\cline{2-15}
		& 0.9$|$0.794 & 0.034 & 20$|$123 & 20125 & 21605 & $-$ & 11 & 30:00 & 23:44 & $-$ & 4.7e-08 & 9.0e-06 & 1.0e-07 & $-$ \\
		\hline
		\multirow{2}{*}{\rotatebox{0}{\Longunderstack{Bosilegrad \\\( (42^{\circ}30'N, 22^{\circ}30'E) \) \\ $(563,73577),\  10511$}}} & 0.0$|$0.316 & 0.026 & 17$|$82 & 20639 & 6783 & 36771 & 37 & 30:00 & 06:32 & 30:00 & 3.6e-08 & 2.9e-06 & 1.0e-07 & 2.7e-03 \\
		\cline{2-15}
		& 0.5$|$0.447 & 0.024 & 23$|$129 & 20267 & 9157 & $-$ & 18 & 30:00 & 09:01 & $-$ & 4.1e-09 & 4.2e-06 & 1.0e-07 & $-$ \\
		\hline
	\end{tabular}
	\caption{Square-root sparse group Lasso model on climate data sets with CV over \eqref{CV2} and \eqref{CV1}. ``$-$'' denotes that the method is not applicable for the instance.}
	\label{tab:climate-4-solvers}
\end{table}

\noindent{\bf Gene data} The first data we use is the {\textit{colon cancer}} data \citep{alon1999broad}\footnote{It is available at \url{http://www.weizmann.ac.il/mcb/UriAlon/download/downloadable-data}.}. This data has been used in \citep{li2017grouped} for the adaptive sparse group Lasso model. It includes 62 tissues (40 colon tumor tissues and 22 normal tissues), and each tissue \(X_i\) includes the expression profiles of 2000 genes. The response variable is assigned according to the label of the tissue: \(Y_i = 1 \) if \(X_i\) is a colon tumor tissue and \(Y_i=-1\) otherwise.

In addition, we use the lung cancer data \citep{monti2003consensus}\footnote{It is available at \url{http://portals.broadinstitute.org/cgi-bin/cancer/publications/view/87}.}, which has also been used in \citep{li2017grouped}. This data includes 197 tissues, and each tissue \(X_i\) includes the expression profiles of 1000 genes. Moreover, 197 tissues are divided into four classes: 17 normal, 139 {{lung adenocarcinoma}}, 21 {{squamous cell carcinomas}}, and 20 {{carcinoids}}.
Corresponding to the latter three classes, we construct three data sets by letting the binary response variables be labels. Specifically, \textit{lung adenocarcinoma} data set has response variable \(Y_i = 1 \) if \(X_i\) is a lung adenocarcinoma tissue and \(Y_i=-1\) otherwise. Similarly, we construct the  \textit{squamous cell carcinomas}  and \textit{carcinoids} data sets.

Another data tested is the acute leukemia data \citep{golub1999molecular}\footnote{It is available at \url{https://github.com/wangyanyanwangyanyan/wangyanyan}.}. This data includes 72 samples, and each sample \(X_i\) includes the expression profiles of 10713 (repeated) genes. Each sample belongs to one of the three classes: BALL, TALL, or AML. The data sets  {\textit{BALL}},  {\textit{TALL}}, and {\textit{AML}} are constructed in the same way as in the lung cancer data above.

In total, we have 7 data sets with binary response variables. For the first 4 data sets,  we used the weighted gene co-expression networks \citep{langfelder2008wgcna} to cluster gene expressions in different groups (modules), and an R package is available for this clustering. The group structures of the last 3 data sets, generated by a similar method, are provided in \citep{li2018grouped}. More details about the group structures can be found in Appendix~\ref{sec:append-gene}.

Table~\ref{tab:gene-4-solvers} presents the results on real gene data sets, where the parameters $w_1$, $w_2$, and $\lambda$ are selected by CV over \eqref{CV1}. For an approximate solution $\beta$, the  classification accuracy is computed by $\left(1-{{\rm nnz} \left( {\rm sign} \left(X_{\rm test}\beta\right) - Y_{\rm test} \right)  }/{N_{\rm test}}\right) \cdot 100\%$.  For problems with \(w_1=0\) where S-TISP is applicable, we also include S-TISP in the table.
The low MSE and high classification accuracy  in Table~\ref{tab:gene-4-solvers} suggest that the square-root sparse group Lasso model is  effective in selecting (groups of) genes related to certain types of tissues and reliable in predicting the class of a tissue;  see Appendix~\ref{sec:append-gene} for more details. Moreover, the selected  \(w_1\not\in\{0,1\} \) by CV \eqref{CV1} in 4 out of 7 data sets indicates  that the {sparse} group Lasso regularizer can improve the performances of classification and gene selection, compared with the Lasso or group Lasso regularizer. For the \textit{BALL} and \textit{TALL} data sets, pADMM and dADMM return overfitting solutions. Besides, one can observe that the running time for different methods is quite similar, since the gene data sets are of small to medium size and they are pre-processed  such that highly irrelevant genes were screened out.

\begin{table}[htbp] \centering
	\fontsize{6pt}{10pt}\setlength{\tabcolsep}{2pt}
	\selectfont
	\renewcommand{\arraystretch}{2}
	\begin{tabular}{|c|r|rr|rrrr|rrrr|rrrr|}
		\hline
		\multicolumn{1}{|c|}{\Longunderstack{problem}} & \multicolumn{1}{c|}{\(w_1|\lambda\)} &  \multicolumn{1}{c}{MSE} &
		\multicolumn{1}{c|}{\Longunderstack{accur}} & \multicolumn{4}{c|}{\Longunderstack{iter }} & \multicolumn{4}{c|}{\Longunderstack{time }} & \multicolumn{4}{c|}{\Longunderstack{error }} \\
		\cline{5-16}
		\multicolumn{1}{|c|}{\Longunderstack{\((N_{\rm train},n), g\) }} & \multicolumn{1}{c|}{ } &  \multicolumn{1}{c}{ } &
		\multicolumn{1}{c|}{\Longunderstack{ }} &
		PP & pA & dA & ST & PP & pA & dA & ST & PP & pA & dA & ST  \\
		\hline
		\multirow{1}{*}{\rotatebox{0}{\Longunderstack{\textit{colon cancer} \\ $(42,2000), 11$ }}} & 0.8$|$0.631 & 0.448 & 90\% & 9$|$13 & 2676 & 665 & $-$ & 00 & 00 & 00 & $-$ & 4.9e-08 & 1.7e-08 & 9.5e-08 & $-$ \\
		\hline
		\multirow{1}{*}{\rotatebox{0}{\Longunderstack{\textit{lung adenocarcinoma} \\ $(132,1000), 8$}}} & 0.0$|$0.562 & 0.156 & 97\% & 9$|$13 & 717 & 613 & 5347 & 00 & 00 & 00 & 00 & 5.7e-09 & 1.4e-08 & 9.9e-08 & 1.0e-07 \\
		\hline
		\multirow{1}{*}{\rotatebox{0}{\Longunderstack{\textit{squamous cell carcinomas} \\ $(131,1000), 8$}}} & 0.0$|$1.259 & 0.131 & 97\% & 10$|$15 & 1861 & 690 & 2783 & 00 & 00 & 00 & 00 & 5.4e-09 & 3.7e-08 & 9.6e-08 & 1.0e-07 \\
		\hline
		\multirow{1}{*}{\rotatebox{0}{\Longunderstack{\textit{carcinoids} \\ $(131,1000), 8$}}} & 0.9$|$1.000 & 0.008 & 100\% & 9$|$12 & 1076 & 509 & $-$ & 00 & 00 & 00 & $-$ & 4.9e-09 & 2.6e-08 & 9.6e-08 & $-$ \\
		\hline
		\multirow{1}{*}{\rotatebox{0}{\Longunderstack{\textit{BALL} \\ $(48,10713), 42$}}} & 0.0$|$0.200 & 0.141 & 100\% & 10$|$19 & 453 & 404 & 1596 & 00 & 00 & 00 & 02 & 1.6e-10 & $*$9.9e-08 & $*$4.9e-08 & $\#$2.1e-13 \\
		\hline
		\multirow{1}{*}{\rotatebox{0}{\Longunderstack{\textit{TALL} \\ $(48,10713), 42$}}} & 0.9$|$0.891 & 0.126 & 96\% & 11$|$25 & 1400 & 2342 & $-$ & 01 & 02 & 03 & $-$ & 3.9e-11 & $*$9.1e-08 & $*$1.0e-07 & $-$ \\
		\hline
		\multirow{1}{*}{\rotatebox{0}{\Longunderstack{\textit{AML} \\ $(48,10713), 42$}}} & 0.9$|$1.585 & 0.193 & 100\% & 12$|$37 & 18393 & 1755 & $-$ & 00 & 29 & 02 & $-$ & 2.4e-08 & 2.1e-08 & 5.7e-08 & $-$ \\
		\hline
	\end{tabular}
	
	\caption{Square-root sparse group Lasso model on gene data sets with CV over \eqref{CV1}. ``$-$'' denotes that the method is not applicable for the instance.
{The errors reported are $\Delta_{\rm kkt}$, $*\Delta_{\rm pd.gap}$, and $\#\Delta_{\rm var.gap} $, given by \eqref{eq:error-kkt}, \eqref{eq:error-gap}, and \eqref{eq:error-change}, respectively.}
}\label{tab:gene-4-solvers}
\end{table}

\subsection{Comparison of Efficiency for Solving the Square-root Fused Lasso Problem}
In this section, we compare PPDNA, pADMM, and  dADMM for solving the square-root fused Lasso problem when the regularizer $p$ is the fused Lasso regularizer \eqref{fused-reg}. We found that a framework of ADMM was applied in \citep{MR4159660} for solving the square-root fused Lasso problem. However, it seems that they applied ADMM in a non-rigorous way as their formulation of the augmented Lagrangian function might not be correct. Except \citep{MR4159660}, there is currently no solver for solving the square-root fused Lasso problem, to the best of our knowledge. For all tables in this section, we denote  PPDNA, pADMM, and dADMM  by ``PP'', ``pA'', and ``dA'', respectively.

\subsubsection{UCI Data}
Again, we use UCI data sets  \textit{housing}, \textit{bodyfat}, \textit{pyrim}, and \textit{triazines} described in Section~\ref{sec:synthetic-group}. We choose \(w_1=w_2=0.5 \)  and \(\lambda \in \left\{1,5,\lambda_{\text{Jia}} \right\} \). See Appendix~\ref{sec:theo-lambda} for the details of \(\lambda_{\text{Jia}}\) given in \citep{MR4159660}.  The numerical results of PPDNA, pADMM, and dADMM are reported in Table~\ref{tab:UCI-fused}. Table~\ref{tab:UCI-fused} shows that our PPDNA substantially outperforms both pADMM and dADMM for solving the square-root fused Lasso problem on the UCI data sets. In particular, our PPDNA {takes} less than 1 minute for all instances; while pADMM {fails} to return accurate solutions within 30 minutes for more than half of the instances. One can conclude that our algorithm is efficient for solving the square-root fused Lasso problem on the UCI data sets.

\begin{table}[p] \centering
	\fontsize{6pt}{10pt}\setlength{\tabcolsep}{3pt} \selectfont
	\begin{tabular}{|r|r|r|rrr|rrr|rrr|}
		\hline
		\multicolumn{1}{|c|}{\Longunderstack{problem  }} & \multicolumn{1}{c|}{\(\lambda\)} & \multicolumn{1}{c|}{\Longunderstack{nnz$|$}} & \multicolumn{3}{c|}{\Longunderstack{iter }} & \multicolumn{3}{c|}{\Longunderstack{time }} & \multicolumn{3}{c|}{\Longunderstack{error }} \\
		\cline{4-12}
		\multicolumn{1}{|c|}{\Longunderstack{\((N,n)\) }} & \multicolumn{1}{c|}{} & \multicolumn{1}{c|}{nnzB} &
		PP & pA & dA & PP & pA & dA & PP & pA & dA  \\
		\hline
		\multirow{3}{*}{\rotatebox{0}{\Longunderstack{\textit{housing} \\ $(253,77520)$}}} & 1.000 & 336$|$129 & 19$|$119 & 35001 & 13873 & 05 & 25:01 & 05:48 & 5.2e-08 & 8.4e-08 & 1.0e-07 \\
		& 5.000 & 168$|$28 & 16$|$81 & 42814 & 30001 & 03 & 30:00 & 11:40 & 2.6e-08 & 2.1e-06 & 5.1e-08 \\
		& 9.282 & 141$|$18 & 14$|$53 & 42015 & 32501 & 02 & 30:00 & 12:22 & 4.0e-08 & 4.7e-07 & 7.5e-08 \\
		\hline
		\multirow{3}{*}{\rotatebox{0}{\Longunderstack{\textit{bodyfat} \\ $(126,116280)$}}} & 1.000 & 295$|$63 & 24$|$121 & 27501 & 6125 & 05 & 15:25 & 02:05 & 2.2e-08 & 9.5e-08 & 1.0e-07 \\
		& 5.000 & 132$|$12 & 14$|$56 & 45001 & 27501 & 02 & 25:36 & 08:34 & 4.3e-09 & 9.1e-08 & 6.8e-08 \\
		& 8.969 & 79$|$6 & 13$|$41 & 54052 & 35001 & 01 & 30:00 & 10:29 & 2.3e-09 & 3.8e-07 & 8.1e-08 \\
		\hline
		\multirow{3}{*}{\rotatebox{0}{\Longunderstack{\textit{pyrim} \\ $(37,169911)$}}} & 1.000 & 535$|$69 & 40$|$169 & 67611 & 116842 & 08 & 30:00 & 30:00 & 4.3e-07 & $*$8.5e-05 & 7.8e-01 \\
		& 5.000 & 386$|$14 & 14$|$52 & 67844 & 20001 & 02 & 30:00 & 05:10 & 4.5e-08 & 8.8e-02 & 3.5e-08 \\
		& 7.848 & 669$|$10 & 15$|$69 & 60001 & 24828 & 03 & 26:31 & 06:17 & 6.1e-10 & 2.6e-08 & 2.5e-08 \\
		\hline
		\multirow{3}{*}{\rotatebox{0}{\Longunderstack{\textit{triazines} \\ $(93,557845)$}}} & 1.000 & 2771$|$165 & 34$|$247 & 8176 & 15422 & 43 & 30:00 & 30:00 & 5.8e-08 & 1.3e-05 & 3.1e-03 \\
		& 5.000 & 1310$|$57 & 22$|$117 & 8122 & 15399 & 22 & 30:00 & 30:00 & 8.9e-08 & 8.6e-05 & 1.2e-04 \\
		& 9.278 & 196$|$12 & 13$|$37 & 8053 & 15716 & 10 & 30:00 & 30:00 & 1.9e-08 & 5.7e-04 & 8.2e-04 \\
		\hline
	\end{tabular}	
	\caption{Square-root fused Lasso model on UCI data sets with \(w_1=w_2=0.5 \).
{The errors reported are $\Delta_{\rm kkt}$ and $*\Delta_{\rm pd.gap}$, given by \eqref{eq:error-kkt} and \eqref{eq:error-gap}, respectively.}
}\label{tab:UCI-fused}
\end{table}

\subsubsection{Real Data}
In this section, we test the square-root fused Lasso model on four real data sets  used in \citep{MR4159660}.
The \textit{inbred mouse}\footnote{It is available at \url{https://www.ncbi.nlm.nih.gov/geo/query/acc.cgi?acc=GSE3330}.} data includes 60 samples, where each sample includes 22689 genes, and the response variable is the feature measured by stearoyl-coenzyme desaturase 1 with probe set ID given as 1415965$\_$at.
The \textit{rat eye}\footnote{It is available at \url{https://www.ncbi.nlm.nih.gov/geo/query/acc.cgi?acc=GSE5680}.} data includes 120 rats samples, where each sample includes 31098 gene probes, and the response variable is selected with respect to 1389163$\_$at.
The \textit{credit card}\footnote{It is available at \url{https://www.kaggle.com/mlg-ulb/creditcardfraud}.} data includes 284807 transactions samples in which 492 samples are labeled as frauds, and each sample {includes} 29 features.
The \textit{safe driver}\footnote{It is available at \url{https://www.kaggle.com/c/porto-seguro-safe-driver-prediction/data}.} data includes 595212 car insurance observations, where each observation includes 57 features.
We refer the readers to \citep[Section~4.3]{MR4159660} for detailed descriptions of the data.

Again, we choose \(w_1=w_2=0.5 \)  and \(\lambda \in \left\{1,5,\lambda_{\text{Jia}} \right\} \) with \(\lambda_{\text{Jia}}\) given in \citep{MR4159660}. The comparisons of PPDNA, pADMM, and dADMM for solving the square-root fused Lasso model on the four real data sets  are reported in Table~\ref{tab:real-fused}.  On can observe that our PPDNA significantly outperforms both pADMM and dADMM for all instances. In particular, both pADMM and dADMM {fail} to solve the problem \textit{safe driver} within 30 minutes.

\begin{table}[p] \centering
	\fontsize{6pt}{10pt}\setlength{\tabcolsep}{3pt} \selectfont
	\begin{tabular}{|r|r|r|rrr|rrr|rrr|}
		\hline
		\multicolumn{1}{|c|}{\Longunderstack{problem  }} & \multicolumn{1}{c|}{\(\lambda\)} & \multicolumn{1}{c|}{\Longunderstack{nnz$|$}} & \multicolumn{3}{c|}{\Longunderstack{iter }} & \multicolumn{3}{c|}{\Longunderstack{time }} & \multicolumn{3}{c|}{\Longunderstack{error }} \\
		\cline{4-12}
		\multicolumn{1}{|c|}{\Longunderstack{\((N,n)\) }} & \multicolumn{1}{c|}{} & \multicolumn{1}{c|}{nnzB} &
		PP & pA & dA & PP & pA & dA & PP & pA & dA  \\
		\hline
		\multirow{3}{*}{\rotatebox{0}{\Longunderstack{\textit{inbred mouse} \\ $(30,22689)$}}} & 1.000 & 85$|$4 & 11$|$21 & 2231 & 6626 & 02 & 05 & 12 & 2.9e-08 & 1.1e-08 & 8.8e-08 \\
		& 5.000 & 1347$|$1 & 13$|$42 & 2043 & 779 & 00 & 04 & 01 & 1.2e-08 & 5.8e-10 & 8.6e-08 \\
		& 6.913 & 3092$|$1 & 13$|$42 & 2009 & 601 & 00 & 05 & 01 & 2.4e-09 & 1.1e-09 & 7.4e-08 \\
		\hline
		\multirow{3}{*}{\rotatebox{0}{\Longunderstack{\textit{rat eye} \\ $(60,31098)$}}} & 1.000 & 433$|$10 & 12$|$30 & 1861 & 2200 & 01 & 09 & 09 & 1.3e-08 & 2.2e-08 & 1.0e-07 \\
		& 5.000 & 3303$|$4 & 12$|$40 & 1549 & 1305 & 00 & 09 & 05 & 4.4e-08 & 3.4e-08 & 1.0e-07 \\
		& 7.798 & 7297$|$4 & 12$|$35 & 1549 & 666 & 00 & 08 & 02 & 2.6e-08 & 7.7e-08 & 1.0e-07 \\
		\hline
		\multirow{3}{*}{\rotatebox{0}{\Longunderstack{\textit{credit card} \\ $(142404,29)$}}} & 1.000 & 27$|$26 & 11$|$24 & 5485 & 1343 & 02 & 42 & 21 & 3.5e-09 & 7.3e-08 & 8.1e-08 \\
		& 5.000 & 25$|$24 & 16$|$84 & 15001 & 5210 & 05 & 01:49 & 01:19 & 6.1e-08 & 3.2e-08 & 6.9e-08 \\
		& 5.683 & 25$|$24 & 18$|$96 & 15001 & 9953 & 06 & 01:52 & 02:28 & 3.2e-08 & 2.7e-08 & 7.6e-08 \\
		\hline
		\multirow{3}{*}{\rotatebox{0}{\Longunderstack{\textit{safe driver} \\ $(595212,57)$}}} & 1.000 & 48$|$35 & 31$|$223 & 35850 & 14429 & 01:29 & 30:00 & 30:00 & 2.7e-08 & 1.9e-05 & 1.3e-04 \\
		& 5.000 & 48$|$30 & 23$|$134 & 36352 & 14333 & 39 & 30:00 & 30:00 & 2.8e-08 & 5.3e-07 & 5.4e-05 \\
		& 6.242 & 47$|$28 & 24$|$133 & 35173 & 14308 & 37 & 30:00 & 30:00 & 1.9e-08 & 4.6e-07 & 1.0e-05 \\
		\hline
	\end{tabular}
	
	\caption{Square-root fused Lasso model on real data sets with \(w_1=w_2=0.5\). }
	\label{tab:real-fused}
\end{table}

\subsection{Simulations with Varying Noise Levels}

\begin{figure}[p]
	\centering
	\includegraphics[width=0.48\textwidth]{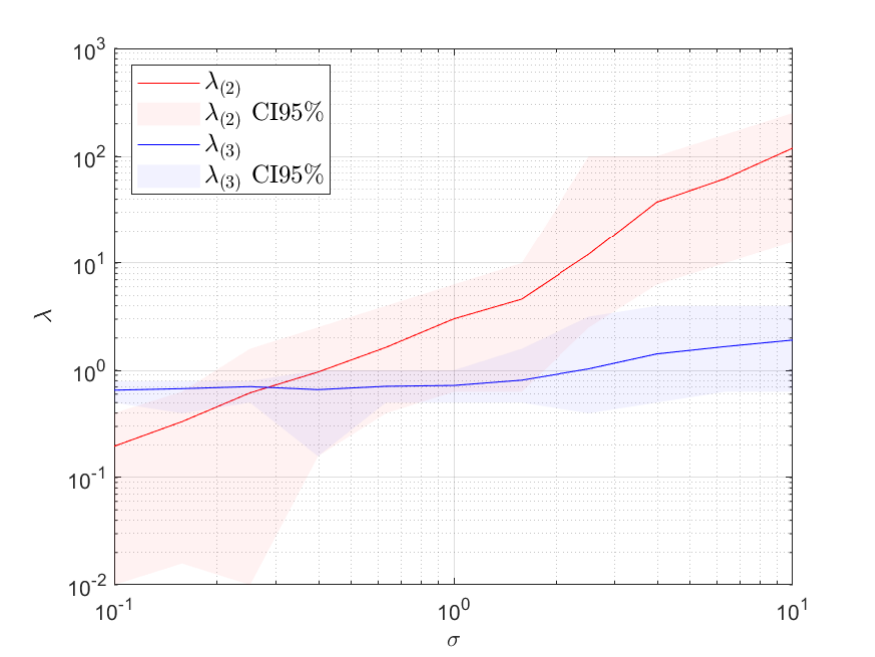}
    \includegraphics[width=0.48\textwidth]{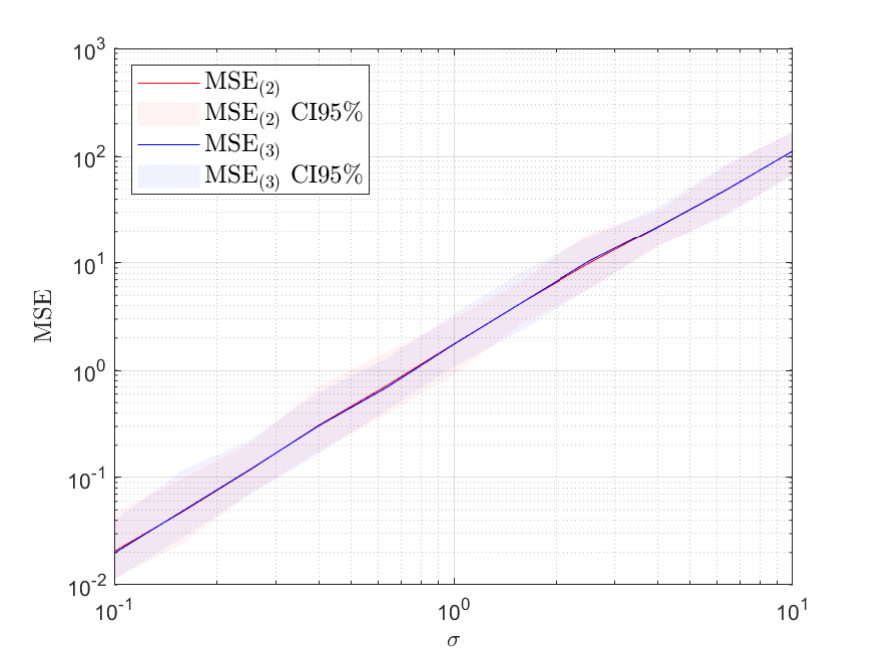}
	\caption{The mean of (left)  \(\lambda\) and (right) MSE, and the $95\%$ confidence interval (CI).}
	\label{fig:robust}
\end{figure}

Here we are interested in how the noise level \(\sigma\) in \eqref{linear-reg} is associated with the tuning parameter \(\lambda\) in the least-square model \eqref{eq:lassotype} and the square-root model \eqref{eq:square-root}. In particular, we let \(p\) be the sparse group Lasso regularizer \eqref{group-reg}. In  this experiment, the noise level \(\sigma\) is chosen from the set \( \left\{10^{-1},10^{-0.8},\dots,10^{1} \right\} \),  the  data is simulated
from \textbf{Example~3} in Section~\ref{sec:synthetic-group} with \(N = 100, g= 50, n = 150 \),
and the parameters $w_1$ and $w_2$ are set to be $0.5$. We regard the 8-fold CV selected tuning parameters \(\lambda_{\eqref{eq:lassotype}}\) and \(\lambda_{\eqref{eq:square-root}}\) as the optimal parameters for model \eqref{eq:lassotype} and model \eqref{eq:square-root}, respectively. We repeat the experiments 100 times, and we plot $\lambda$ and MSE against $\sigma$ in Figure~\ref{fig:robust}.

We can observe from the left panel of  Figure~\ref{fig:robust} that as $\sigma$ varies from \(10^{-1}\) to \(10^{1}\), the curve of \(\lambda_{\eqref{eq:square-root}}\) roughly remains flat with the values of \(\lambda_{\eqref{eq:square-root}}\) staying in the small range  \([10^{-0.5},10^{0.5}]\). In contrast, \(\lambda_{\eqref{eq:lassotype}}\) shown in the left panel of  Figure~\ref{fig:robust} varies wildly from \(10^{-1}\) to \(10^{2}\). This phenomenon verifies numerically the advantage of the square-root model \eqref{eq:square-root} in terms of choosing the tuning parameters compared with the least-square model \eqref{eq:lassotype}, namely, the parameter $\lambda$ in the square-root model \eqref{eq:square-root} can be tuned independent of the noise level $\sigma$.
In addition, it can be observed from {the slope of the MSE curves in }the right panel of Figure~\ref{fig:robust} that both the MSE achieved by \eqref{eq:lassotype} and that by \eqref{eq:square-root}  are {approximately quadratic} in the noise level $\sigma$.

\section{Conclusion} \label{sec:conclusion}
In this paper, we have given a unified proof to show that any square-root regularized model whose penalty function being the sum of a simple norm and a seminorm can be interpreted as the  distributionally robust optimization formulation of the corresponding  least-squares  problem. For solving a generic square-root regularized model, we have developed a proximal point dual semismooth Newton algorithmic framework
to efficiently solve the resulting convex minimization problem whose objective is the sum of two nonsmooth terms
corresponding to the square-root loss and the regularizer respectively. We have illustrated that the general framework can be adopted to solve the square-root sparse group Lasso and the square-root fused Lasso models.  Our extensive numerical experiments have shown that the proposed algorithm is indeed highly efficient for solving the square-root sparse group Lasso and the square-root fused Lasso models, as compared to popular first order methods based on the ADMM framework.


\acks{The research of Kim-Chuan Toh is supported by the Ministry of Education, Singapore, under its Academic Research Fund Tier 3 grant call (MOE-2019-T3-1-010).
The research of Yangjing Zhang is supported by the National Natural Science Foundation of China under grant number 12201617.
}



\appendix
\section{Remark on Proposition~\ref{prop:seminorm}}\label{sec:remark}

Proposition~\ref{prop:seminorm} includes a similar result in \citep{maurer2012structured}.
For simplicity, we show that a simplified form of their results can be derived from Proposition~\ref{prop:seminorm}. In \citep{maurer2012structured}, $\mathcal{M}:= \{M_{(l)}\}_{1\leq l \leq g}$ denotes a set of symmetric matrices $M_{(l)} \in\mathbb{S}^n$, and the operator $\|\cdot\|_{\mathcal{M}}\,:\mathbb{R}^n\to [0,+\infty]$ is defined by
\begin{equation} \label{eq:M-norm}
	\|\alpha\|_{\mathcal{M}} := \inf_{{\alpha}_{(1)},\dots,{\alpha}_{(g)}}
\left\{ \sum_{l=1}^{g} \|{\alpha}_{(l)}\| \,\bigg|\, {\alpha}_{(l)} \in\mathbb{R}^n, \, \sum_{l=1}^{g}M_{(l)}{\alpha}_{(l)}=\alpha \right\}\,\,\forall \,\alpha \in \mathbb{R}^n.
\end{equation}	
It was shown in \citep{maurer2012structured} that $\|\cdot\|_{\mathcal{M}}$ is indeed a norm on the subspace of $\mathbb{R}^n$ where it is finite, and the dual norm is given by
\begin{equation}\label{eq:M-dual-norm}
\|\beta\|_{\mathcal{M}^*} := \sup_{1\leq l \leq g} \left\{ \|M_{(l)}\beta\|  \right\}\,\,\forall \,\beta \in \mathbb{R}^n.
\end{equation}
We define $B:= \left(M_{(1)},\dots,M_{(g)}\right)^T \in\mathbb{R}^{ng\times n}$ and $Q(\bar{\alpha}) := \sup_{1\leq l \leq g} \left\{ \|\alpha_{(l)}\|  \right\}$ for any $\bar{\alpha} = (\alpha_{(1)},\dots,\alpha_{(g)})\in\mathbb{R}^n \times \cdots \times \mathbb{R}^n$. $Q$ is a norm on $\mathbb{R}^{ng}$, and  its dual norm is given by $Q_*(\bar{\alpha}) := \sum_{l=1}^{g} \|\alpha_{(l)}\|\,\,\forall\,\bar{\alpha} = (\alpha_{(1)},\dots,\alpha_{(g)})\in\mathbb{R}^n \times \cdots \times \mathbb{R}^n$. We can see from \eqref{eq:M-norm} that $\|\cdot\|_{\mathcal{M}}$ coincides with $p_*$ in \eqref{eq:pstar2}, i.e.,
$
\|\alpha\|_{\mathcal{M}} = \inf_{\bar{\alpha}} \{Q_*(\bar{\alpha})\,|\,B^T \bar{\alpha} = \alpha\}\,\,\forall \,\alpha \in \mathbb{R}^n.
$
It is indeed a norm on $\mathtt{Range}(B^T)$ from Proposition~\ref{prop:seminorm}(a).
We can also see from \eqref{eq:M-dual-norm} that $\|\cdot\|_{\mathcal{M}^*}$ coincides with $p$ in \eqref{eq:pform2}, i.e.,
$
\|\beta\|_{\mathcal{M}^*} = Q(B\beta)\,\,\forall \,\beta \in \mathbb{R}^n.
$
Therefore, we can derive the result that $\|\cdot\|_{\mathcal{M}}$ and $\|\cdot\|_{\mathcal{M}^*}$ are dual to each other directly from Proposition~\ref{prop:seminorm}.

\section{Theoretical Parameter $\lambda$} \label{sec:theo-lambda}

The independence of the tuning parameter $\lambda$ on the unknown noise level $\sigma$ is one of the nice statistical properties of the square-root regularized model. In this section, we present various selections of $\lambda$ which has been studied in literature for the problem \eqref{eq:square-root}. The cumulative distribution function of the standard normal distribution $\mathcal{N}(0,1)$ is denoted by $\Phi(x):=\frac{1}{\sqrt{2\pi}} \int_{-\infty}^x e^{-t^2/2} dt$.
The cumulative distribution function of the $F$-distribution with the degrees of freedom $a$ and $b$ is denoted by $\mathcal{F}_{a,b}$.
The quantile function with respect to a cumulative distribution function $F$ is denoted by $\mathcal{Q}_F(a) := \inf \{x\in\mathbb{R} \mid a \leq F(x) \} $.

\textbf{Sparse group Lasso regularizer \eqref{group-reg}}
$p(\beta) = w \|\beta\|_1 + (1-w)\sum_{l=1}^g \sqrt{\abs{G_l}} \|\beta_{G_l}\|\,\,\forall\,\beta\in\mathbb{R}^n,$
where \(w\in[0,1]\). We summarize the selections of $\lambda$ which are independent on $\sigma$ in \citep{belloni2011square,bunea2013group,stucky2017sharp,blanchet2019robust,blanchet2017distributionally}. In the following formulations, $a$ is chosen to be 0.05. When $N<10^4$, we calculate ${\rm E}[X^TX] \;{\approx}\; \frac{1}{N}X^TX $; otherwise, we randomly sample $10^4$ predictor vectors from \(X\) to form \(\tilde{X} \), and estimate  ${\rm E}[X^TX] \;{\approx}\; \frac{1}{10^4}\tilde{X}^T\tilde{X} $.
\begin{enumerate}
	\item When $w=1 $, $\lambda$ can be selected from $\Lambda_S := \{\lambda_{\text{Bel}},\lambda_{\text{StS}},\lambda_{\text{BlS}} \} $.
	\begin{itemize}
		\item[-] \citep{belloni2011square}  $\lambda_{\text{Bel}} := 1.1\Phi^{-1}\left( 1-\frac{a}{2n} \right)$.
		\item[-] \citep{stucky2017sharp} Denote $t := \sqrt{\log\left(\frac{4}{a}\right)},
		\Delta := \sqrt{ 1 - t \sqrt{\frac{4}{N}}}$,  then $ \lambda_{\text{StS}}:= \sqrt{2}\frac{t}{\Delta}  + \sqrt{2}\left(2+\sqrt{\log(n)}\right) $.
		\item[-] \citep{blanchet2019robust} Estimate   $Z\sim \mathcal{N}(0,{\rm E}[X^TX])   $ and  $F$ as the cumulative distribution of $\frac{\pi}{\pi-2}\norm{Z}_{\infty}^2 $. Let $\hat{\eta}_{1-a}:= \mathcal{Q}_F(1-a) $, and $\lambda_{\text{BlS}}:=\sqrt{\hat{\eta}_{1-a}}$.
	\end{itemize}
	\item When $w=0$,  $\lambda$ can be selected from $\Lambda_G := \{\lambda_{\text{Bun}},\lambda_{\text{StG}},\lambda_{\text{BlG}} \} $.
	\begin{itemize}
		\item[-] \citep{bunea2013group} Denote $T_{\max}:= \max_{1\leq l \leq g}\{\abs{G_l}\},T_{\min}:= \min_{1\leq l \leq g}\{\abs{G_l}\}$,  $\tau_0 := \mathcal{F}^{-1}_{T_{\min},N-T_{\min}}\left( 1-\frac{a}{g}\right)$ and $ \zeta_{\max} := \max \left\{ \frac{1}{N} \norm{X_{G_l } }^2 \mid 1\leq l \leq g \right\}$. Suppose that $T_{\min}\tau_0 + N -T_{\max} >0$, then $\lambda_{\text{Bun}} :=  \sqrt{\frac{\zeta_{\max}\tau_0 }{ T_{\min}\tau_0 + N -T_{\max} }} \sqrt{N}$.
		\item[-] \citep{stucky2017sharp} Denote $t := \sqrt{\log\left(\frac{4}{a}\right)},
		\Delta := \sqrt{ 1 - t \sqrt{\frac{4}{N}}}$,  then $ \lambda_{\text{StG}}:= \sqrt{2}\frac{t}{\Delta}  + \sqrt{2}\left(2+\sqrt{\log(g)}\right) $.
		\item[-] \citep{blanchet2017distributionally} Estimate   $Z\sim \mathcal{N}(0,{\rm E}[X^TX])   $ and  $F$ as the cumulative distribution of $\frac{\pi}{\pi-2}(p_*(Z))^2 $. Let $\hat{\eta}_{1-a}:= \mathcal{Q}_F(1-a) $, and $\lambda_{\text{BlG}}:=\sqrt{\hat{\eta}_{1-a}}$.
	\end{itemize}
\end{enumerate}

\textbf{Fused Lasso regularizer \eqref{fused-reg}}
Following \citep{MR4159660}, we let \(a = 0.05,\,t = \sqrt{{4\log(\frac{1}{a})}/{N}}+{4\log(\frac{1}{a})}/{N}\) and calculate $\lambda_{\text{Jia}} = 2.2\sqrt{2\log(n)/(1+t)}$.

\section{Gene Data Sets} \label{sec:append-gene}

In this section, we provide more details about the gene data sets used in Section~\ref{sec:real-data-group} and their numerical results.

\begin{table}[htbp]
\centering
	\begin{tabular}{|ccccccccccc|}
		\hline
		 black &	blue&	brown	&green&	grey&	magenta&	pink&	purple&	red&	turquoise&	yellow\\
		\hline
		71	&420&	367&	98&	18&	63&	64&	43&	87&	507&	262	\\
		\hline
	\end{tabular}
\caption{The number of genes in each module in the colon cancer data.}
\label{table:colcan}
\end{table}
\begin{figure}[htbp]
	\centering
	\includegraphics[width=0.3\textwidth]{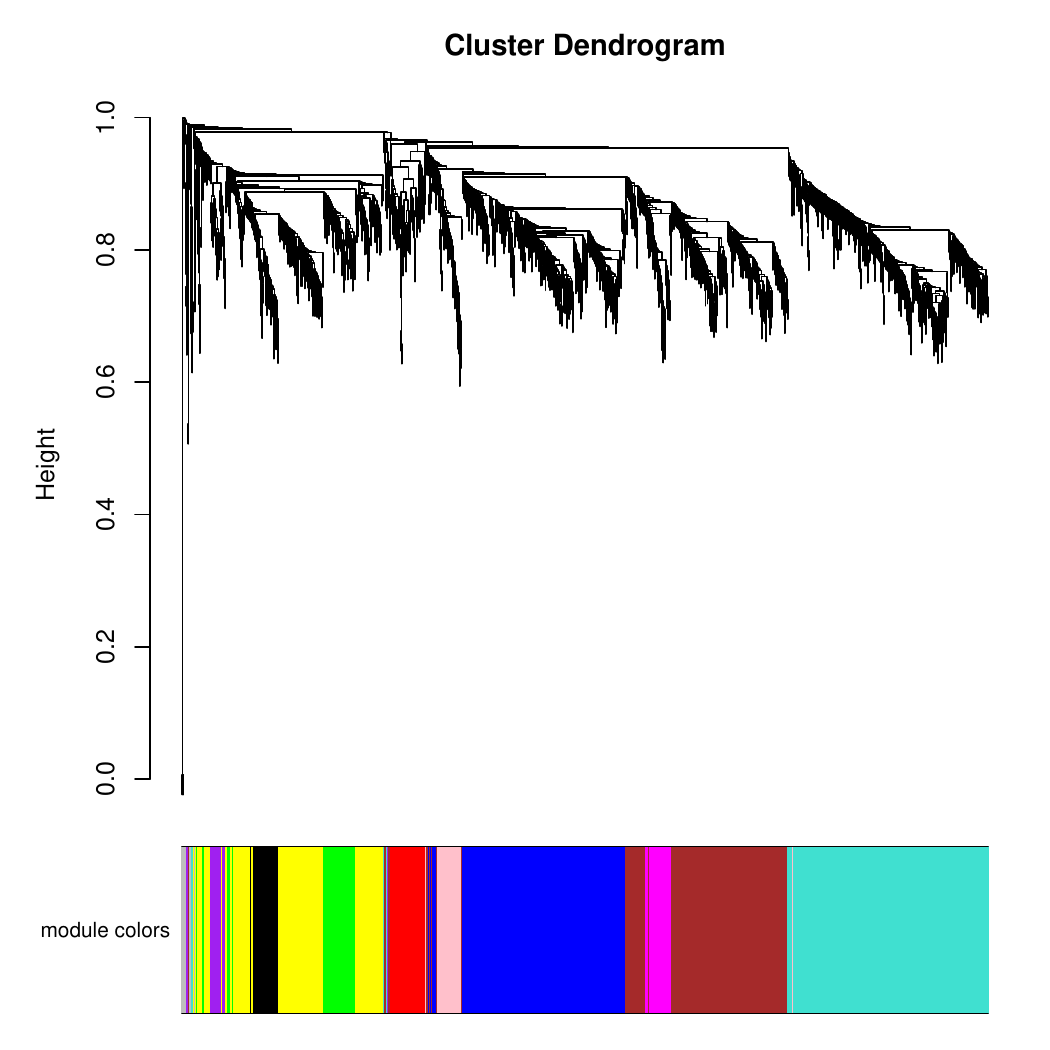}
	\includegraphics[width=0.3\textwidth]{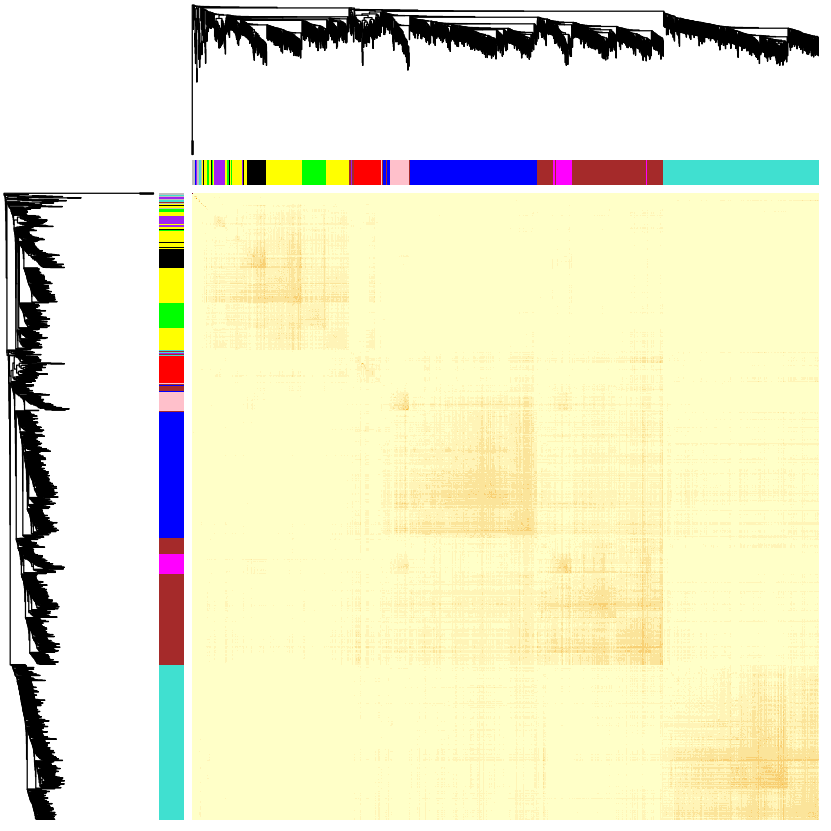}
	\caption{Left: clustering dendrogram of 2000 gene profiles from the colon cancer patients. Right: heatmap of gene-gene connectivity.}
	\label{fig:colon-cancer-data}
\end{figure}
\begin{figure}[htbp]
	\centering
	\includegraphics[width=0.6\textwidth]{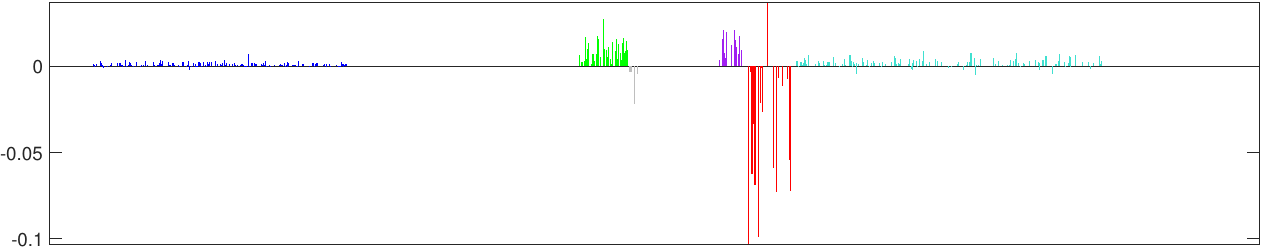}
	\caption{The value of \(\beta\) (coordinates reordered by its groups) in colon cancer data.}
	\label{fig:colon-cancer-beta}
\end{figure}
\begin{table}[htbp] \centering
	\begin{tabular}{|cccccccc|}
		\hline
		black &     blue &    brown &    green  &    grey    &   red &turquoise   & yellow \\\hline
		28     &  186  &     112    &    58    &   172     &   38     &  332   &     74  \\\hline
	\end{tabular}
\caption{The number of genes in each module in the lung cancer data.}
\label{table:luncan}
\end{table}
\begin{figure}[htbp]
	\centering
	\includegraphics[width=0.3\textwidth]{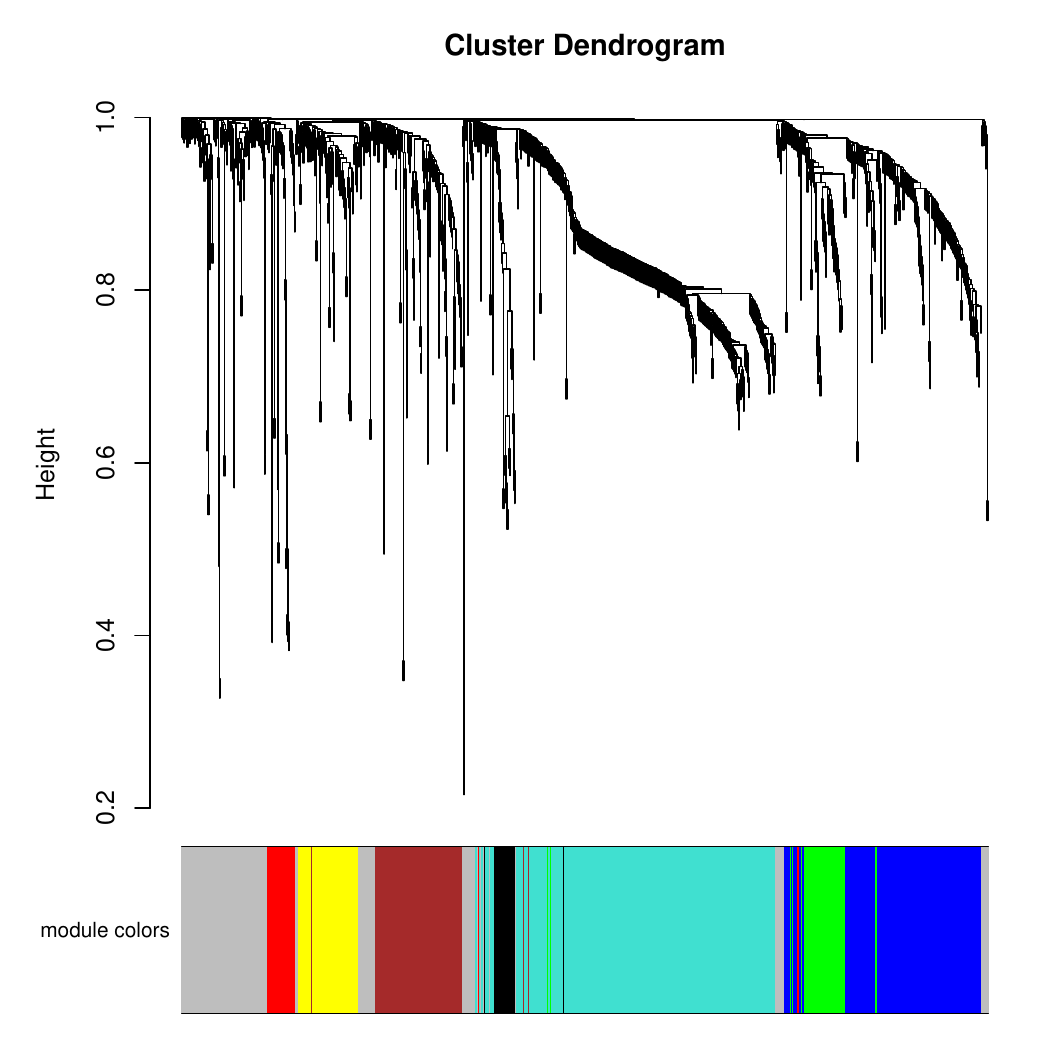}
	\includegraphics[width=0.3\textwidth]{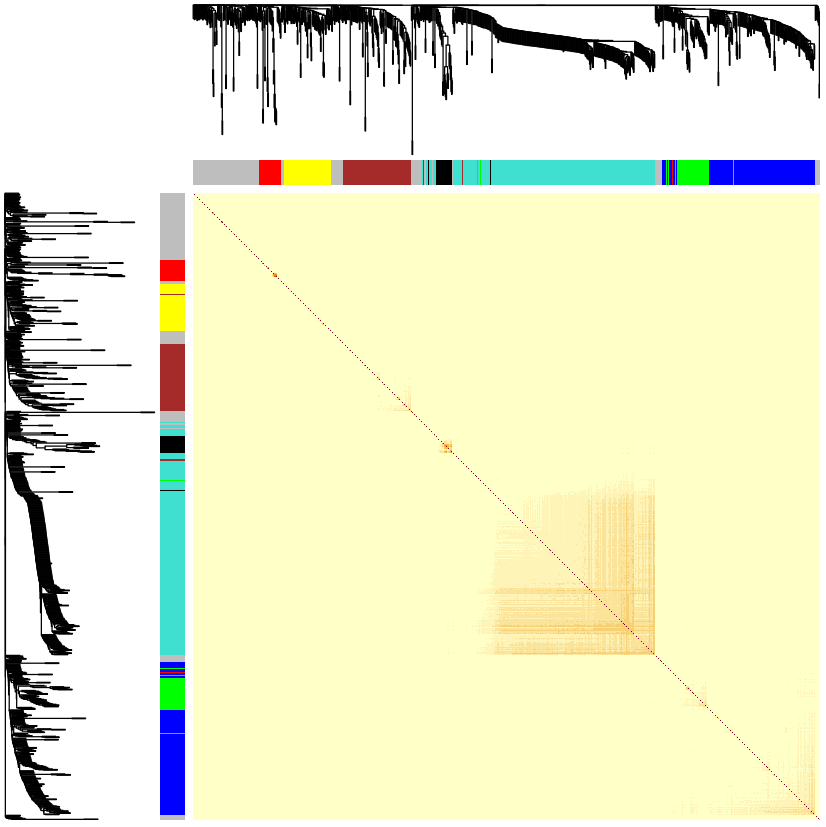}
	\caption{Left: clustering dendrogram of 1000 gene profiles from the lung cancer patients. Right: heatmap of gene-gene connectivity.}
	\label{fig:lung-cancer-data}
\end{figure}
\begin{figure}[htbp]
	\centering
	\includegraphics[width=0.6\textwidth]{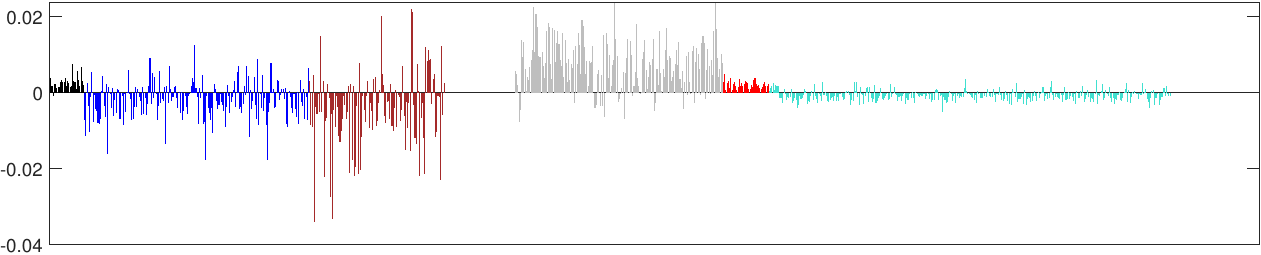}
	\includegraphics[width=0.6\textwidth]{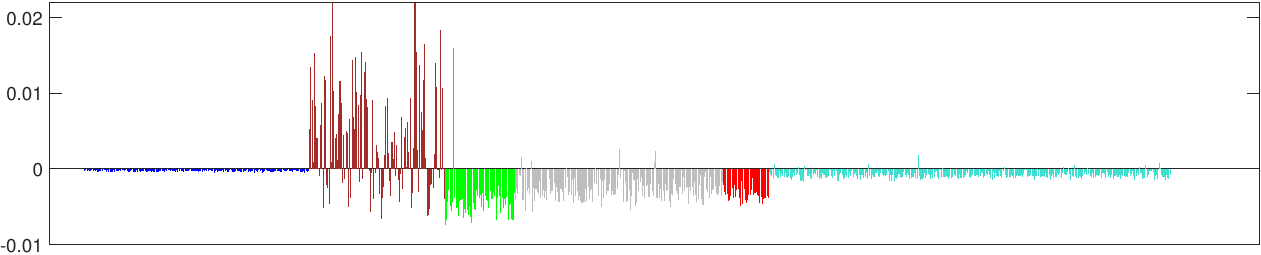}
	\includegraphics[width=0.6\textwidth]{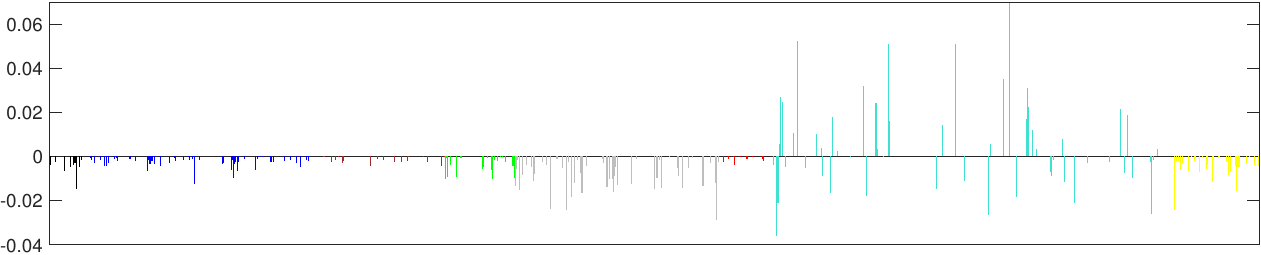}
	\caption{The value of \(\beta\) (coordinates reordered by its groups) in lung cancer data. Top: {\it lung adenocarcinoma}, middle: {\it squamous cell carcinomas}, bottom: {\it carcinoids}.}
	\label{fig:lung-cancer-beta}
\end{figure}
\begin{figure}[h]
	\centering
	\includegraphics[width=0.6\textwidth]{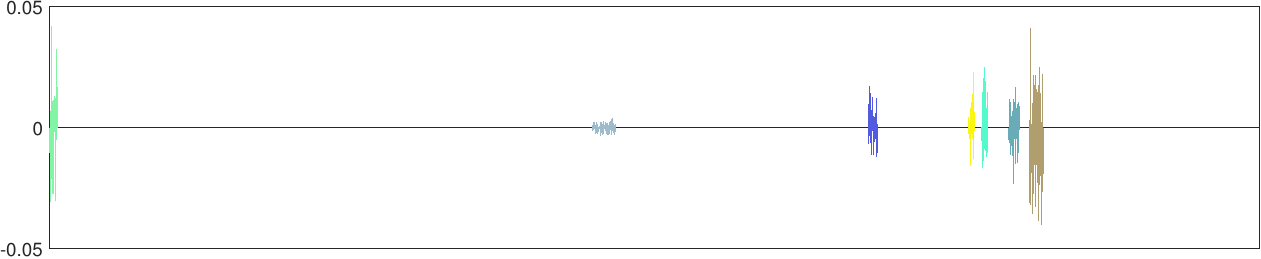}
	\includegraphics[width=0.6\textwidth]{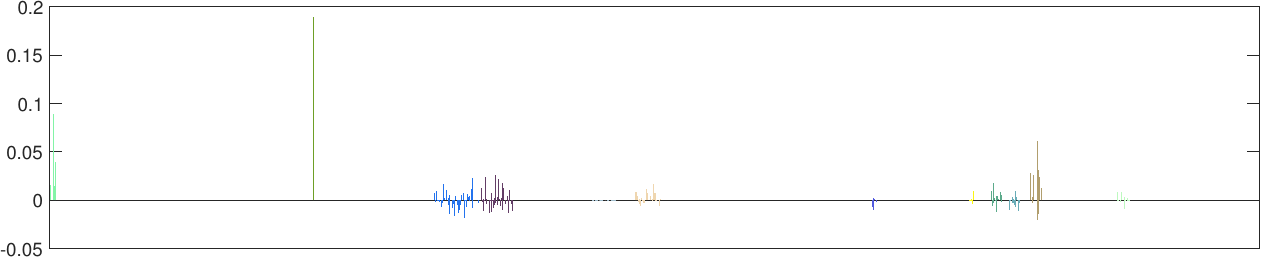}
	\includegraphics[width=0.6\textwidth]{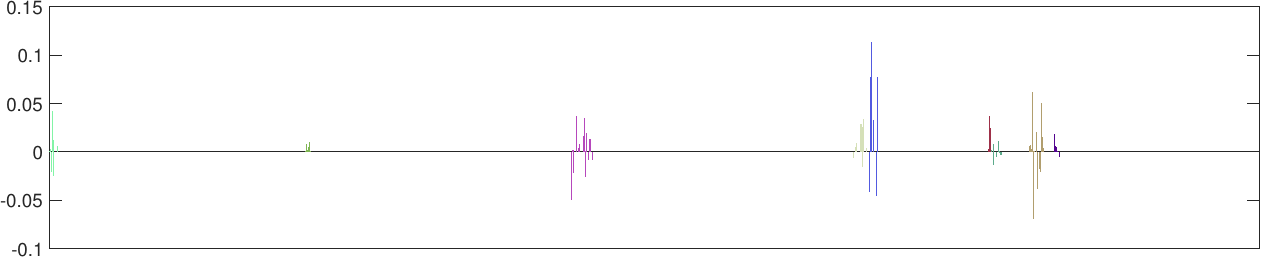}
	\caption{The value of \(\beta\) (coordinates reordered by its groups)  for acute leukemia data. Top: {\it BALL}, middle: {\it TALL}, bottom: {\it AML}.}
	\label{fig:acute-leukemia-beta}
\end{figure}

\textbf{\textit{Colon cancer data}}
This data includes 62 tissues with binary labels, and each tissue includes the expression profiles of 2000 genes. We use the R package WGCNA \citep{langfelder2008wgcna} with the power parameter 6 to construct the weighted gene co-expression networks, and divide 2000 genes into 11 modules which are marked by 11 colors in Figure~\ref{fig:colon-cancer-data}. The number of genes in each module is given in Table~\ref{table:colcan}. The value of the approximate solution \(\beta\) (coordinates reordered by its groups) with CV \eqref{CV1} selected parameters  is presented in Figure~\ref{fig:colon-cancer-beta}. One can observe that there are 6 active groups out of total 11 groups, which may imply that the red, green, and purple groups are key groups to classify the colon cancer patients. Moreover, we can see that \(\beta\) is also sparse within each group, which indicates the effectiveness of the square-root sparse group Lasso regularizer.

\textit{\textbf{Lung cancer data}}
This data includes 197 tissues, and each tissue  includes the expression profiles of 1000 genes. We use the R package WGCNA \citep{langfelder2008wgcna} with the power parameter 6 to construct the weighted gene co-expression networks, and divide 1000 genes into 8 modules which are marked by 8 colors in Figure~\ref{fig:lung-cancer-data}. The number of genes in each module is given in Table~\ref{table:luncan}. The value of the approximate solution \(\beta\) (coordinates reordered by its groups) with CV \eqref{CV1} selected parameters  in \textit{lung adenocarcinoma}, \textit{squamous cell carcinomas}, and \textit{carcinoids} data sets is presented in Figure~\ref{fig:lung-cancer-beta}. One can observe that the set of selected groups varies for different data sets. For example, the brown group is one of the key modules to classify lung adenocarcinoma and squamous cell carcinomas tissues, while the turquoise module is preferred in the carcinoids data set. Note that \eqref{CV1} returns \(w_1=0 \) in the first two data sets and  \(w_1=0.9\) in the last data set, which results in the sparsity of the last regression vector.

\textbf{\textit{Acute leukemia data}}
The raw acute leukemia data set includes 72 samples of 3571 gene expressions. Following \citep{li2018grouped}, the grouping strategy with repeated genes was applied to get the data with 10713 gene expressions.  The value of the approximate solution \(\beta\) (coordinates reordered by its groups) with CV \eqref{CV1} selected parameters  in \textit{BALL}, \textit{TALL}, and \textit{AML} data sets is presented in Figure~\ref{fig:acute-leukemia-beta}. Similar conclusions could be made as in the lung cancer data above.

\newpage
\bibliography{ref}

\end{document}